\documentclass[11pt,a4paper,reqno]{amsart}

\usepackage{amssymb}
\usepackage{mathtools}
\usepackage{amsmath}
\usepackage{bbm}

\usepackage{graphicx}
\usepackage{tikz}
\usetikzlibrary{decorations.pathmorphing}

\usepackage{enumitem}

\usepackage{cite}
\usepackage{etoolbox}
\usepackage{xcolor}
\usepackage{microtype}

\usepackage{hyperref}
\usepackage{prettyref}

\numberwithin{equation}{section}

\definecolor{LinkBlue}{RGB}{45,110,170}

\hypersetup{
  unicode=true,
  colorlinks=true,
  linkcolor=LinkBlue,
  citecolor=LinkBlue,
  urlcolor=LinkBlue
}

\makeatletter
\renewcommand{\subsection}{%
  \@startsection{subsection}{2}%
    {\z@}%
    {2ex plus 1ex minus .2ex}%
    {1ex plus .2ex}%
    {\normalfont\bfseries}%
}
\makeatother

\newcommand{\theoremheadfont}{%
  \fontencoding{T1}%
  \fontfamily{ntxtlf}%
  \fontseries{b}%
  \fontshape{sc}%
  \selectfont
}

\newtheoremstyle{smallcaps}
  {10pt plus 2pt minus 2pt}       
  {10pt plus 2pt minus 2pt}       
  {\itshape}      
  {}              
  {\theoremheadfont} 
  {.}             
  {5pt plus 1pt minus 1pt} 
  {}              

\newtheoremstyle{smallcaps-roman}
  {10pt plus 2pt minus 2pt}
  {10pt plus 2pt minus 2pt}
  {\normalfont}
  {}
  {\theoremheadfont}
  {.}
  {5pt plus 1pt minus 1pt}
  {}

\theoremstyle{smallcaps}

\newtheorem{theorem}{Theorem}[section]
\newtheorem{proposition}[theorem]{Proposition}
\newtheorem{lemma}[theorem]{Lemma}
\newtheorem{corollary}[theorem]{Corollary}

\theoremstyle{smallcaps-roman}

\newtheorem{definition}[theorem]{Definition}
\newtheorem{remark}[theorem]{Remark}
\newtheorem*{assumption*}{Assumption}

\setlist[enumerate]{
  leftmargin=0pt,
  itemindent=2.5em,
  labelsep=0.5em,
  topsep=9pt,
  itemsep=9pt,
  parsep=0pt
}

\usepackage{prettyref}

\newrefformat{theorem}{{\normalfont\textsc{Theorem}~\ref*{#1}}}
\newrefformat{proposition}{{\normalfont\textsc{Proposition}}~\ref*{#1}}
\newrefformat{lemma}{{\normalfont\textsc{Lemma}}~\ref*{#1}}
\newrefformat{corollary}{\textsc{Corollary}~\ref*{#1}}
\newrefformat{definition}{\textsc{Definition}~\ref*{#1}}
\newrefformat{remark}{\textsc{Remark}~\ref*{#1}}

\newcommand{\assumptionref}{%
  \hyperref[assumption:coeff]{{\normalfont\textsc{Assumption}}}%
}

\newcommand{\refsp}[1]{\hyperref[#1]{\prettyref{#1}}}

\newcommand{\thm}{theorem}
\newcommand{\prop}{proposition}
\newcommand{\lem}{lemma}
\newcommand{\cor}{corollary}
\newcommand{\defi}{definition}

\newcommand{\rem}{remark}

\makeatletter

\DeclareRobustCommand\widecheck[1]{%
  {\mathpalette\@widecheck{#1}}%
}

\def\@widecheck#1#2{%
  \setbox\z@\hbox{\m@th$#1#2$}%
  \setbox\tw@\hbox{%
    \m@th
    $#1%
    \widehat{%
      \vrule\@width\z@\@height\ht\z@
      \vrule\@height\z@\@width\wd\z@
    }$%
  }%
  \dp\tw@-\ht\z@
  \@tempdima\ht\z@
  \advance\@tempdima2\ht\tw@
  \divide\@tempdima\thr@@
  \setbox\tw@\hbox{%
    \raise\@tempdima\hbox{%
      \scalebox{1}[-1]{%
        \lower\@tempdima\box\tw@
      }%
    }%
  }%
  {%
    \ooalign{%
      \box\tw@\cr
      \box\z@
    }%
  }%
}

\makeatother

\newcommand{\norm}[1]{%
  \left\lVert #1 \right\rVert
}

\newcommand{\lnorm}[1]{%
  \left\lVert #1 \right\rVert_{L^2}
}

\newcommand{\lpnorm}[2]{%
  \left\lVert #2 \right\rVert_{L^{#1}}
}

\newcommand{\hsnorm}[1]{%
  \left\lVert #1 \right\rVert_{H^s}
}

\newcommand{\lloc}[1]{%
  L^{#1}_{\mathrm{loc}}
}

\newcommand{\ctext}[1]{%
  \raise0.2ex\hbox{%
    \textcircled{\scriptsize #1}%
  }%
}

\newcommand{\aeae}{\mathrm{a.e.}}

\newcommand{\supp}[1]{%
  \mathrm{supp}\,\,#1\,\,
}

\newcommand{\esssup}[1]{%
  \underset{#1}{\mathrm{ess\,sup}}\,
}

\newcommand{\set}[2]{%
  \left\{\,#1\,\middle|\,#2\,\right\}
}

\newcommand{\jb}[1]{%
  \left\langle #1 \right\rangle
}

\newcommand{\Z}{\mathbb{Z}}

\newcommand{\R}{\mathbb{R}}
\newcommand{\C}{\mathbb{C}}

\newcommand{\cA}{\mathcal{A}}
\newcommand{\cB}{\mathcal{B}}

\newcommand{\F}{\mathcal{F}}
\newcommand{\cK}{\mathcal{K}}

\newcommand{\Sc}{\mathcal{S}}

\newcommand{\step}{{S\textsc{tep}}}

\newlist{stepenum}{enumerate}{1}

\setlist[stepenum]{%
  label=\step\ \arabic*.,
  ref=\arabic*
}

\usepackage{fancyhdr}



\allowdisplaybreaks[4]

\title[LWP for generalized INLS]
{Local well-posedness for generalized inhomogeneous
nonlinear Schrödinger equations}

\author{Kohei Inada}

\address{
Graduate School of Advanced Science and Engineering,
Hiroshima University, \\
1-3-1 Kagamiyama, Higashi-Hiroshima, 739-8526 Japan
}

\email{{kohei-inada}@hiroshima-u.ac.jp}

\makeatletter
\patchcmd{\@setauthors}
  {\MakeUppercase{\authors}}
  {\MakeUppercase{\authors}\par
   \vspace{0.5em}
   \begingroup
     \def\author##1{}
     \def\address##1##2{{\normalfont\scriptsize ##2\par}}
     \def\curraddr##1##2{}
     \def\email##1##2{{\normalfont\scriptsize\ttfamily ##2\par}}
     \def\urladdr##1##2{{\normalfont\scriptsize\ttfamily ##2\par}}
     \addresses
   \endgroup}
  {}{}

\let\@setaddresses\@empty
\makeatother



\begin{document}

\begin{abstract}
  We consider the initial value problem for the nonlinear Schrödinger equation
  \begin{equation*}
    i\partial_t u+\Delta u=a(x)|u|^\alpha u,
  \end{equation*}
  where \(a(x)\) is a spatially dependent coefficient. We study a general class of coefficients that includes the usual inhomogeneous nonlinearity \(a(x)=\lambda|x|^{-b}\) and coefficients with multiple singularities. We establish local well-posedness in \(H^s(\mathbb{R}^N)\) for both the \(H^s\)-subcritical and \(H^s\)-critical cases under suitable assumptions on \(a(x)\). We also extend a fractional difference estimate of Fujiwara to Lorentz spaces.
\end{abstract}

\maketitle


\section{Introduction} \label{sec:intro}

\subsection{Background} \label{sec:1-1}

In this paper, we consider the initial value problem for the nonlinear Schr\"odinger equation with a spatially dependent nonlinear term,
\begin{equation} \label{eq:1-1}
  i\partial_tu + \Delta u = a(x)|u|^{\alpha}u, \quad u(0,\cdot) = u_0.
\end{equation}
Here, $u_0\in H^s(\R^N)$, $u:\R\times\R^N\to\C$ is a complex-valued function, and $a:\R^N\to\C$ is a spatially dependent coefficient. We also denote by $\Delta$ the Laplace operator and assume that $s\geq 0$, $N\geq 1$, and $\alpha>0$. 

When $a(x)\equiv \lambda$, equation \refsp{eq:1-1} becomes the classical nonlinear Schr\"o\-dinger equation (NLS) with a power-type nonlinearity,
\begin{equation} \label{eq:NLS}
  i\partial_tu + \Delta u = \lambda|u|^{\alpha}u, \quad u(0,\cdot) = u_0. \tag{NLS}
\end{equation}
The NLS arises naturally in nonlinear optics and serves, in particular, as a model for the propagation of high-intensity laser beams in a homogeneous medium with a Kerr nonlinearity. It is also one of the prototypical nonlinear dispersive equations and has been extensively studied from a mathematical point of view.

On the other hand, when $a(x)=\lambda|x|^{-b}$, equation \refsp{eq:1-1} becomes the inhomogeneous nonlinear Schr\"odinger equation (INLS),
\begin{equation} \label{eq:INLS}
  i\partial_tu + \Delta u = \lambda|x|^{-b}|u|^{\alpha}u, \quad u(0,\cdot) = u_0. \tag{INLS}
\end{equation}
Beginning in the late twentieth century, the guiding of high-power laser beams through a plasma channel formed by a preceding laser pulse has been studied; see, for example, Gill~\cite{gil}. More directly related to spatially inhomogeneous nonlinearities, Fibich--Wang~\cite{F-W} considered the nonlinear Schr\"odinger equation with a spatially dependent coefficient $V(x)$,
\begin{equation} \label{eq:1-2}
  i\partial_tu + \Delta u + V(x)|u|^{\alpha}u = 0.
\end{equation}
In de Bouard--Fukuizumi~\cite{D-F}, $V(x)$ is interpreted as a quantity proportional to the electron density in plasma, and equation \refsp{eq:1-2} is studied under the assumption $V(x)\sim |x|^{-b}\ (x\to\infty)$. Moreover, Genoud--Stuart~\cite{G-S} studied standing wave solutions for general $V(x)$. They also showed that the power-law nonlinearity $|x|^{-b}|u|^\alpha u$ naturally appears through a suitable scaling limit.

Genoud--Stuart~\cite{G-S} also discussed the $H^1$ well-posedness of \refsp{eq:INLS}. More precisely, using the abstract argument of Cazenave~\cite{Caz}, they obtained local well-posedness roughly in the range $0<b<\min\{2,N\}$, $0<\alpha<\frac{4-2b}{N-2}$. To the best of our knowledge, this was the first well-posedness result for \refsp{eq:INLS}. Later, Guzmán~\cite{Guz} developed a Strichartz-based well-posedness theory in $H^s$, $0\leq s\leq1$, under suitable restrictions on $b$ and $\alpha$. The critical case $\alpha = \frac{4-2b}{N-2s}$ was first treated by Kim--Lee--Seo~\cite{K-L-S}. For $N\geq3$, they obtained well-posedness results by using weighted $L^p$ spaces under some restrictions on $\alpha$ and $b$. In a different approach, Aloui--Tayachi~\cite{A-T} introduced a framework based on Lorentz spaces, using the fact that $|x|^{-b}\in L^{\frac{N}{b},\infty}$. Their method treats the subcritical and critical cases in a unified way. They also established results on unconditional uniqueness and continuous dependence on the initial data.

More recently, Campos--Correia--Farah~\cite{C-C-F} considered a wider range than \cite{A-T} and established local well-posedness in the range 
\begin{equation*}
  0<b<\min\left\{2,N-s, \frac{N}{2} + 1 - s\right\}.
\end{equation*}
They also examined whether the assumptions on $b$ and $\alpha$, which arise from technical requirements, are essential. In particular, they showed that a certain form of ill-posedness occurs outside the corresponding well-posedness range, under additional assumptions.

Very recently, several results concerning this equation have appeared. Lee~\cite{lee} studied the three-dimensional energy-critical INLS with strong singularity and established local and small-data global well-posedness for $3/2\leq b<11/6$ by using improved inhomogeneous Strichartz estimates in weighted spaces. Bhimani--Dhingra--Sohani~\cite{B-D-S} investigated the INLS in Fourier--Lebesgue and modulation spaces and obtained low-regularity well-posedness results in these settings. Grira--Miyazaki--Tayachi~\cite{G-M-T} considered inhomogeneous nonlinear Schrödinger equations without gauge invariance and studied the threshold for the existence of scattering states by using non-admissible Strichartz estimates in Lorentz spaces.

{ 
The aim of this study is to extend the single power-type coefficient $a(x)=\lambda|x|^{-b}$, which has mainly been considered in previous studies on INLS, to a more general class of coefficients $a(x)$, and to establish local well-posedness in the subcritical and critical cases. By imposing assumptions on $a(x)$ and its fractional derivative in Lorentz spaces, we can treat more general coefficients, including bounded terms and several terms with different singularities. Thus, our results establish local well-posedness in a unified framework not only for the standard INLS but also for the classical NLS and equations with several power-type singularities.

In addition, we prove local well-posedness by a different method based on Lorentz spaces, under the same conditions on the singularity and nonlinearity exponents as those obtained by Campos--Correia--Farah \cite{C-C-F}. Lorentz spaces serve as a natural framework for handling the singularities of power-type coefficients. Their properties allow us to estimate the nonlinear terms without relying on the explicit form of the coefficient. By establishing the key estimates for the nonlinear terms in the framework of Lorentz spaces, we not only extend the previous results for INLS but also unify the proofs for both subcritical and critical cases.}

\subsection{Main Results} \label{sec:1-2}

As explained in the previous subsection, we study the local well-posedness of equation \refsp{eq:1-1} with a general spatially dependent coefficient $a(x)$, not necessarily of the power form $|x|^{-b}$. In this subsection, we first discuss the subcritical and critical conditions that are closely related to the main results, and then state the main theorems. 

For fixed $N\geq 1$ and $\alpha>0$, we define
\begin{equation} \label{eq:1-3}
  s_c(b) \coloneq \frac{N}{2} - \frac{2 - b}{\alpha}, \quad 0\leq b<2.
\end{equation}
In the typical case \refsp{eq:INLS} with $a(x) = \lambda|x|^{-b}$, consider the scaling $u_{\kappa}(t,x) = \kappa^{\frac{2-b}{\alpha}}u(\kappa^2 t, \kappa x)$. Then
\begin{equation*}
  \norm{u_{\kappa}(0)}_{\dot{H}^s} = \kappa^{\frac{2-b}{\alpha} + s - \frac{N}{2}}\norm{u_0}_{\dot{H}^s}.
\end{equation*}
Thus, the $\dot{H}^s$-norm is invariant under this scaling when $s = s_c(b)$. Accordingly, the equation is called $H^s$-subcritical when $s>s_c(b)$, $H^s$-critical when $s=s_c(b)$, and $H^s$-supercritical when $s<s_c(b)$. In the range considered in the main theorems, $b<2$, $s<\frac{N}{2}$, $\alpha>0$, a direct calculation gives
\begin{gather*}
  s>s_c(b) \iff 0<\alpha<\frac{4-2b}{N-2s}, \\
  s=s_c(b) \iff \alpha = \frac{4-2b}{N-2s}, 
\end{gather*}
and 
\begin{equation*}
  s<s_c(b) \iff \alpha > \frac{4-2b}{N-2s}.
\end{equation*}
For this reason, our main results are stated in terms of conditions on $\alpha$: \refsp{\thm:1-1} treats the subcritical case, while \refsp{\thm:1-2} deals with the critical case. 

For a general coefficient $a(x)$, the equation \refsp{eq:1-1} does not necessarily have the scaling invariance described above. Following the assumption in \refsp{sec:2-4}, we fix an exponent $b \geq 0$ and set $s_c\coloneq s_c(b)$. By \refsp{eq:1-3}, we have $s_c(b_{\ast})\leq s_c$ for every $b_{\ast}\in\{0,b\}$. 

Before stating the main results, we define the admissible pairs used in this paper. For $2\leq q,r\leq \infty$, a pair $(q,r)$ is called admissible if 
\begin{equation*}
    \frac{2}{q} + \frac{N}{r} = \frac{N}{2}
\end{equation*} 
and $(q,r,N)\neq (2,\infty, 2)$. We also define
\begin{equation*}
    \cA_0 \coloneq \set{(q,r):\text{admissible pair}}{2<r<\frac{2N}{N-2}}
\end{equation*}
where, for $N=1,2$, the condition on $r$ is replaced by $2<r<\infty$. The more general notion of an $\dot{H}^s$-admissible pair and the corresponding class $\cA_{s}$ are introduced later in \refsp{\defi:adm}. 

In the proofs, for each $b_{\ast}\in\{0,b\}$, we will choose an admissible pair $(\gamma_{b_{\ast}},\rho_{b_{\ast}})\in\cA_0$ suited to the corresponding nonlinear estimate. The precise choice is given in \refsp{sec:4} and \refsp{sec:5}. In the critical case with $b>0$, the terms corresponding to $b_{\ast}=0$ are treated as in the subcritical case, while the critical estimates use $b_{\ast} = b$. When $b=0$, there is only one exponent. 

In this paper, local well-posedness in $H^s$ means the existence of a solution on a short time interval, uniqueness, and continuous dependence on the initial data. We now state the main results for the subcritical and critical cases. 

\begin{theorem}[{Local well-posedness in $H^s$ in the subcritical case}] \label{theorem:1-1}
  Let $N\geq 1$, $0\leq s\leq 1$, and $s<\frac{N}{2}$. Suppose that $a(x)$ satisfies the \assumptionref\ in \refsp{sec:2-4} for some $b\geq 0$. Assume further that 
  \begin{equation*}
    0\leq b < \min\left\{ 2, N-s, \frac{N}{2}+1-s \right\}, \quad 0<\alpha < \frac{4-2b}{N-2s}.
  \end{equation*}
  Then, for each $b_{\ast}\in\{0,b\}$, there exists an admissible pair $(\gamma_{b_{\ast}},\rho_{b_{\ast}})\in\cA_0$ such that, for every $u_0\in H^s(\R^N)$, there exist $T_{\max}(u_0)>0$ and a unique solution $u$ to \refsp{eq:1-1} satisfying
  \begin{equation*}
    u\in C([0,T_{\max}(u_0)),H^s(\R^N))\cap \bigcap_{b_{\ast}\in\{0,b\}} \lloc{\gamma_{b_{\ast}}}([0,T_{\max}(u_0));H_{\rho_{b_{\ast}},2}^{s}(\R^N)).
  \end{equation*}
  Moreover, the following properties hold. 
  \begin{enumerate}[label=(\roman*)]
    \item For every $0<T<T_{\max}(u_0)$, the solution $u$ is unique in the above class. 
    \item For every $0<T<T_{\max}(u_0)$ and every admissible pair $(q,r)$ with $r<\infty$, $u\in L^{q}(0,T;H_{r,2}^{s}(\R^N))$. 
    \item For every $R>0$, there exists $T=T(R)>0$ satisfying the following: for every admissible pair $(q,r)\in\cA_0$ and every $\varepsilon>0$, there exists $\delta = \delta(R,\varepsilon,q,r) > 0$ such that, whenever $u_0, v_0\in H^s(\R^N)$ satisfy
    \begin{equation*}
      \hsnorm{u_0}, \hsnorm{v_0}\leq R, \quad \hsnorm{u_0 - v_0} < \delta,
    \end{equation*} 
    the corresponding solutions $u$ and $v$ to \refsp{eq:1-1} satisfy
    \begin{equation} \label{eq:6}
      \sup_{0\leq t\leq T}\norm{u(t)-v(t)}_{H^{s}} + \norm{u-v}_{L^q(0,T;H_{r,2}^{s})} < \varepsilon.
    \end{equation}
  \end{enumerate}
\end{theorem}

We next state the result for the critical case. 

\begin{theorem}[{Local well-posedness in $H^s$ in the critical case}] \label{theorem:1-2}
  Let $N\geq 1$, $0\leq s\leq 1$, and $s<\frac{N}{2}$. Suppose that $a(x)$ satisfies the \assumptionref in \refsp{sec:2-4} for some $b\geq 0$. Assume further that 
  \begin{equation*}
    0\leq b < \min\left\{ 2, N-s, \frac{N}{2}+1-s \right\}, \quad \alpha = \frac{4-2b}{N-2s}.
  \end{equation*}
  Then, for each $b_{\ast}\in\{0,b\}$, there exists an admissible pair $(\gamma_{b_{\ast}},\rho_{b_{\ast}})\in\cA_0$ such that, for every $u_0\in H^s(\R^N)$, there exist $T_{\max}(u_0)>0$ and a unique solution $u$ to \refsp{eq:1-1} satisfying
  \begin{equation*}
    u\in C([0,T_{\max}(u_0)),H^s(\R^N))\cap \bigcap_{b_{\ast}\in\{0,b\}} \lloc{\gamma_{b_{\ast}}}([0,T_{\max}(u_0));H_{\rho_{b_{\ast}},2}^{s}(\R^N)).
  \end{equation*}
  Moreover, the following properties hold. 
  \begin{enumerate}[label=(\roman*)]
    \item For every $0<T<T_{\max}(u_0)$, the solution $u$ is unique in the above class. 
    \item For every $0<T<T_{\max}(u_0)$ and every admissible pair $(q,r)$ with $r<\infty$, $u\in L^{q}(0,T;H_{r,2}^{s}(\R^N))$. 
    \item For every $(q,r)\in\cA_0$, every compact set $\cK\subset H^s(\R^N)$, and every $\varepsilon>0$, there exist $T=T(\cK)>0$ and $\delta = \delta(\cK,\varepsilon,q,r) > 0$ such that, whenever $u_0, v_0\in \cK$ satisfy
    \begin{equation*}
      \hsnorm{u_0 - v_0} < \delta,
    \end{equation*} 
    the corresponding solutions $u$ and $v$ to \refsp{eq:1-1} satisfy
    \begin{equation} \label{eq:cri}
      \sup_{0\leq t\leq T}\norm{u(t)-v(t)}_{H^{s}} + \norm{u-v}_{L^q(0,T;H_{r,2}^{s})} < \varepsilon.
    \end{equation}
  \end{enumerate}
\end{theorem}

In this paper, a solution to \refsp{eq:1-1} is understood as a mild solution satisfying the corresponding Duhamel formula, that is, a solution to the integral equation
\begin{equation} \label{eq:1-5}
  u(t) = \Phi(u)(t) \coloneq e^{it\Delta}u_0 - i\int_0^t e^{i(t-\tau)\Delta}(a(x)|u|^\alpha u(\tau))d\tau.
\end{equation}
 In what follows, we refer to such a mild solution simply as a solution. 

To prove the main theorems, we combine the Strichartz estimates with nonlinear estimates in Lorentz spaces and show that $\Phi$ is a contraction on a suitable function space. The local solution is then constructed by the Banach fixed-point theorem. 

In \refsp{sec:2}, we introduce the notation used throughout the paper and collect basic estimates in Lorentz and Sobolev--Lorentz spaces, together with the Strichartz estimates. We also discuss the class of coefficients $a(x)$ considered in this paper in \refsp{sec:2-4}.

When $0<\alpha<1$, the derivative of the nonlinear mapping may fail to be Lipschitz continuous, and a Hölder estimate is used instead. As a result, a more precise difference estimate is needed to prove continuous dependence on the initial data. For this purpose, in \refsp{sec:3} we extend Fujiwara's fractional derivative estimate~\cite{fuj} for differences to the setting of Lorentz spaces. 

In the subcritical case, the nonlinear estimates contain a small factor coming from the length of the time interval. Thus, the contraction property follows by taking the interval sufficiently small. In the critical case, this small factor is lost for the component corresponding to the critical exponent $b$. Therefore, we use the smallness of the Strichartz norm on a short time interval instead. This difference is the reason why the subcritical and critical cases are treated separately. Moreover, in the proof of continuous dependence in the critical case, this smallness must be obtained uniformly with respect to the initial data. Accordingly, \refsp{\thm:1-2}(iii) is formulated for compact subsets of $H^s$. 

We prove \refsp{\thm:1-1} in \refsp{sec:4}. In \refsp{sec:5}, we prove \refsp{\thm:1-2} by combining the extended Strichartz estimate with the choice of exponents needed in the critical case.

\section{Preliminaries} \label{sec:2}

\subsection{Notation} \label{sec:2-1}

We collect the notation used throughout the paper. The symbol $C$ denotes a positive constant that may vary from line to line. For positive quantities $a$ and $b$, we write $a\lesssim b$ if $a\leq Cb$ for some constant $C>0$ independent of $a$ and $b$. When necessary, the parameters on which the constant may depend are indicated by subscripts. For $1\leq p\leq \infty$, we denote by $p^\prime$ the H\"older conjugate exponent of $p$, that is, $\frac{1}{p} + \frac{1}{p^\prime} = 1$. We define the Fourier transform and its inverse by
\begin{equation*}
    \F[f](\xi) = \widehat{f}(\xi) \coloneq \int_{\R^N} e^{-ix\cdot\xi} f(x) dx,
\end{equation*}
and
\begin{equation*}
    \F^{-1}[f](x) = \widecheck{f}(x) \coloneq \frac{1}{(2\pi)^N} \int_{\R^N} e^{ix\cdot\xi} f(\xi) d\xi.
\end{equation*}
For $s\in\R$, we define the inhomogeneous and homogeneous differentiation operators by
\begin{equation*}
    \jb{D}^sf \coloneq \F^{-1}\left[ \jb{\xi}^s \widehat{f}(\xi) \right], \quad D^sf \coloneq \F^{-1}\left[ |\xi|^s \widehat{f}(\xi) \right],
\end{equation*}
where $\jb{\xi} = (1 + |\xi|^2)^{\frac{1}{2}}$. For $1\leq p < \infty$, the usual Lebesgue space $L^p(\R^N)$ is equipped with the norm
\begin{equation*}
    \lpnorm{p}{f} \coloneq \left( \int_{\R^N} |f(x)|^p dx \right)^{\frac{1}{p}}.
\end{equation*}
For $p=\infty$, we use the usual modification:
\begin{equation*}
    \lpnorm{\infty}{f} \coloneq \esssup{x\in\R^N} |f(x)|.
\end{equation*}
We denote by $\Sc(\R^N)$ the Schwartz space and by $\Sc^{\prime}(\R^N)$ its dual, the space of tempered distributions. The Sobolev space $H^s(\R^N)$ is defined by
\begin{equation*}
    H^s(\R^N) \coloneq \set{f \in \Sc^{\prime}(\R^N)}{\jb{D}^sf \in L^2(\R^N)},
\end{equation*}
with norm $\norm{f}_{H^s} \coloneq \norm{\jb{D}^sf}_{L^2}$. For $s\geq 0$, the homogeneous Sobolev space $\dot{H}^s(\R^N)$ is defined by
\begin{equation*}
    \dot{H}^s(\R^N) \coloneq \set{f \in \Sc^{\prime}(\R^N)}{D^sf \in L^2(\R^N)}
\end{equation*}
with semi-norm $\norm{f}_{\dot{H}^s} \coloneq \norm{D^sf}_{L^2}$. Unless otherwise stated, all function spaces are defined over $\R^N$. Thus, for example, $L^p$, $H^s$, and $\dot{H}^s$ stand for $L^p(\R^N)$, $H^s(\R^N)$, and $\dot{H}^s(\R^N)$, respectively. For a time interval $I\subset\R$ and a Banach space $X$, $L^q(I;X)$ denotes the Lebesgue space of $X$-valued functions with norm
\begin{equation*}
    \norm{f}_{L^q(I;X)} \coloneq \left( \int_I \norm{f(t)}_X^q dt \right)^{\frac{1}{q}}.
\end{equation*}
For $q = \infty$, we use the usual modification. 

\subsection{Lorentz spaces and Sobolev--Lorentz spaces} \label{sec:2-2}

In this subsection, we recall the definitions and basic properties of Lorentz spaces and Sobolev--Lorentz spaces used in this paper. For further details on Lorentz spaces, see, for example, \cite{gra,one}. 

Let $f$ be a measurable function on $\R^N$. We define its distribution function $d_f$ on $(0,\infty)$ by
\begin{equation*}
  d_f(\lambda) = \left|\set{x\in\R^N}{|f(x)|>\lambda}\right|.
\end{equation*}
Here, $|A|$ denotes the Lebesgue measure of a subset $A\subset \R^N$. The rearrangement $f^\ast$ of $f$ is then defined by
\begin{equation*}
  f^\ast(s) = \inf\set{\lambda>0}{d_f(\lambda)\leq s}.
\end{equation*}
Clearly, both $d_f$ and $f^\ast$ are nonnegative decreasing functions. 

\begin{definition}[{Lorentz spaces}]
  Let $p\in (0,\infty)$ and $q\in (0,\infty]$. The Lorentz space $L^{p,q}(\R^N)$ is defined by
  \begin{equation*}
    L^{p,q}(\R^N) = \set{f:\R^N\to\C\ \text{measurable}}{\norm{f}_{L^{p,q}}<\infty},
  \end{equation*}
  where
  \begin{equation*}
    \norm{f}_{L^{p,q}} = \begin{dcases}
      \left( \frac{q}{p} \int_0^{\infty} \left( s^{1/p}f^\ast(s) \right)^q \frac{1}{s} ds \right)^{1/q}, & 0<q<\infty, \\
      \sup_{s>0} s^{1/p}f^\ast(s), & q=\infty.
    \end{dcases}
  \end{equation*}
\end{definition}

The functional $\norm{\ \cdot\ }_{L^{p,q}}$ is in general only a quasi-norm, that is, the triangle inequality holds up to a multiplicative constant. Thus, $L^{p,q}$ is a quasi-normed space, and in fact a quasi-Banach space. 

For fixed $p$, the Lorentz spaces become larger as $q$ increases, that is, $L^{p,q_1}\hookrightarrow L^{p,q_2}$ and there exists a constant $C = C(p,q_1,q_2) > 0$ such that
\begin{equation*} 
  \norm{f}_{L^{p,q_2}} \leq C \norm{f}_{L^{p,q_1}}
\end{equation*}
for $1\leq q_1 \leq q_2 \leq \infty$. In particular, $L^{p,1}\hookrightarrow L^{p}\hookrightarrow L^{p,\infty}$. Moreover, for every $0 < p,\mu < \infty$ and $0<q\leq\infty$, 
\begin{equation*} 
  \norm{|f|^{\mu}}_{L^{p,q}} \leq C \norm{f}_{L^{\mu p,\mu q}}^\mu.
\end{equation*}

For $1<p<\infty$ and $1\leq q\leq\infty$, one can introduce an equivalent norm on $L^{p,q}$. Define
\begin{equation*}
  f^{\ast\ast}(t) = \frac{1}{t}\int_0^t f^\ast(s) ds, \quad t>0,
\end{equation*}
and
\begin{equation*}
  \norm{f}_{L^{p,q}}^\ast = \begin{dcases}
    \left( \frac{q}{p} \int_0^{\infty} \left( s^{1/p}f^{\ast\ast}(s) \right)^q \frac{1}{s} ds \right)^{1/q}, \quad 1\leq q<\infty, \\
    \sup_{s>0} s^{1/p}f^{\ast\ast}(s), \quad q=\infty.
  \end{dcases}
\end{equation*}
Then $L^{p,q}$, equipped with $\norm{f}_{L^{p,q}}^\ast$, is a Banach space, and
\begin{equation*}
  \norm{f}_{L^{p,q}} \leq \norm{f}_{L^{p,q}}^\ast \leq \frac{p}{p-1} \norm{f}_{L^{p,q}}.
\end{equation*}
Hence, the two quantities $\norm{f}_{L^{p,q}}$ and $\norm{f}_{L^{p,q}}^\ast$ are equivalent for $1<p<\infty$ and $1\leq q<\infty$. When $p=q$, $L^{p,p} = L^{p}$, and, with the usual modification, $L^{\infty, \infty} = L^{\infty}$. It is also well known that $L^{p,\infty}$ coincides with the weak $L^p$ space. For example, $|x|^{-b}\notin L^{p}(\R^N)$ for every $p$, while $|x|^{-b}\in L^{\frac{N}{b}, \infty}$. 

In addition, we use the interpolation inequality 
\begin{equation} \label{eq:2-3}
  \norm{f}_{L^{p,q}} \leq C \norm{f}_{L^{p_1,q_1}}^{\theta} \norm{f}_{L^{p_2,q_2}}^{1-\theta}, 
\end{equation}
where $1 < p, p_1, p_2 < \infty$, $1 \leq q, q_1, q_2 \leq \infty$, $\theta\in (0,1)$, and
\begin{equation*}
  \frac{1}{p} = \frac{\theta}{p_1} + \frac{1-\theta}{p_2}, \hspace{20pt} \frac{1}{q} \leq \frac{\theta}{q_1} + \frac{1-\theta}{q_2}.
\end{equation*}

For $s\geq 0$, $1 < p < \infty$, and $1 \leq q \leq \infty$, we define the Sobolev--Lorentz spaces (see \cite[page 571]{H-Y-Z}) by 
\begin{equation*}
  H_{p,q}^s(\R^N) = \set{f\in\Sc^\prime(\R^N)}{\jb{D}^sf\in L^{p,q}}, 
\end{equation*}
and
\begin{equation*}
  \dot{H}_{p,q}^s(\R^N) = \set{f\in\Sc^\prime(\R^N)}{D^sf\in L^{p,q}}.
\end{equation*}
We equip $H_{p,q}^s(\R^N)$ with the norm
\begin{equation*}
  \norm{f}_{H_{p,q}^s} = \norm{f}_{L^{p,q}} + \norm{\jb{D}^sf}_{L^{p,q}},
\end{equation*}
and we equip $\dot{H}_{p,q}^s(\R^N)$ with the semi-norm
\begin{equation*}
  \norm{f}_{\dot{H}_{p,q}^s} = \norm{D^sf}_{L^{p,q}}.
\end{equation*}
In particular, $H_{p,q}^s$ is a reflexive Banach space if $1<p,q<\infty$.

We next collect several inequalities in Lorentz spaces that will be used to estimate the nonlinear terms.

\begin{lemma}[{H\"older's inequality in Lorentz spaces \protect\cite{one}}] \label{lemma:Holder}
  Let $1 < p_1, p_2 \leq \infty$, $1 \leq q_1, q_2 \leq \infty$, and suppose that $f\in L^{p_1,q_1}$, $g\in L^{p_2,q_2}$. Then the following statements hold.
  \begin{enumerate}[label=(\roman*)]
    \item If $\frac{1}{p_1} + \frac{1}{p_2} < 1$, then, for $p$ satisfying 
    \begin{equation*}
      \frac{1}{p} = \frac{1}{p_1} + \frac{1}{p_2}, 
    \end{equation*}
    and any $q\geq 1$ satisfying
    \begin{equation*}
      \frac{1}{q} \leq \frac{1}{q_1} + \frac{1}{q_2},
    \end{equation*}
    we have $fg\in L^{p,q}$. Moreover, there exists a constant 
    \begin{equation*}
      C = C(p,q,p_1,p_2,q_1,q_2) > 0
    \end{equation*}
    such that
    \begin{equation*}
      \norm{fg}_{L^{p,q}} \leq C \norm{f}_{L^{p_1,q_1}} \norm{g}_{L^{p_2,q_2}}.
    \end{equation*}
    \item If $\frac{1}{p_1} + \frac{1}{p_2} = 1$, $\frac{1}{q_1} + \frac{1}{q_2} \geq 1$, then $fg\in L^1$. Furthermore, there exists a constant $C = C(p_1,p_2,q_1,q_2) > 0$ such that
    \begin{equation*}
      \norm{fg}_{1} \leq C \norm{f}_{L^{p_1,q_1}} \norm{g}_{L^{p_2,q_2}}.
    \end{equation*}
    \item If $p_2 = q_2 = \infty$, then $fg\in L^{p_1,q_1}$, and 
    \begin{equation*}
      \norm{fg}_{L^{p_1,q_1}} \leq C \norm{f}_{L^{p_1,q_1}} \norm{g}_{\infty}
    \end{equation*}
    for some constant $C > 0$. 
  \end{enumerate}
\end{lemma}

\begin{lemma}[{Fractional Leibniz rule in Lorentz spaces\protect\cite[Theorem 6.1]{C-N}}] \label{lemma:Leibniz}
  Let $s\in [0,1]$, and suppose that $1 < p, p_1, p_2, p_3, p_4 < \infty$, $1 \leq q, q_1, q_2, q_3, q_4 \leq \infty$ satisfy
  \begin{equation*}
    \frac{1}{p} = \frac{1}{p_1} + \frac{1}{p_2} = \frac{1}{p_3} + \frac{1}{p_4}, 
  \end{equation*}
  and 
  \begin{equation*}
    \frac{1}{q} \leq \frac{1}{q_1} + \frac{1}{q_2} , \quad \frac{1}{q} \leq \frac{1}{q_3} + \frac{1}{q_4}.
  \end{equation*}
  Then
  \begin{equation*}
    \norm{D^s(fg)}_{L^{p,q}} \lesssim \norm{D^sf}_{L^{p_1,q_1}}\norm{g}_{L^{p_2,q_2}} + \norm{f}_{L^{p_3,q_3}}\norm{D^sg}_{L^{p_4,q_4}}.
  \end{equation*}
\end{lemma}

\begin{lemma}[{Fractional chain rule in Lorentz spaces\protect\cite[Lemma 2.4]{A-T}}] \label{lemma:chain}
  Let $s\in [0,1]$ and $F\in C^1(\C,\C)$. Suppose that $1 < p, p_1, p_2 < \infty$, $1 \leq q, q_1, q_2 < \infty$ satisfy
  \begin{equation*}
    \frac{1}{p} = \frac{1}{p_1} + \frac{1}{p_2}, \quad \frac{1}{q} \leq \frac{1}{q_1} + \frac{1}{q_2}.
  \end{equation*}
  Then 
  \begin{equation*}
    \norm{D^sF(f)}_{L^{p,q}} \lesssim \norm{F^\prime(f)}_{L^{p_1,q_1}}\norm{D^sf}_{L^{p_2,q_2}}.
  \end{equation*}
\end{lemma}

\begin{lemma}[{Fractional chain rule for H\"older continuous functions in Lorentz spaces \cite[Proposition 2.7]{par}}] \label{lemma:chain2}
  Let $0<\alpha\leq 1$, and let $F$ be an $\alpha$-H\"older continuous function. Suppose that $1 < p, p_1, p_2 < \infty$, $1 \leq q, q_1, q_2 < \infty$ satisfy
  \begin{equation*}
    0 < \frac{1}{p} = \frac{\alpha - s/s_0}{p_1} + \frac{s/s_0}{p_2} < 1, \quad \text{and} \quad 0 \leq \frac{1}{q} \leq \frac{\alpha - s/s_0}{q_1} + \frac{s/s_0}{q_2} \leq 1.
  \end{equation*}
  Then, for $0<s<\alpha$ and $\frac{s}{\alpha} < s_0 < 1$, we have
  \begin{equation*}
    \norm{D^sF(f)}_{L^{p,q}} \lesssim \norm{f}_{L^{p_1,q_1}}^{\alpha - s/s_0} \norm{D^{s_0}f}_{L^{p_2,q_2}}^{s/s_0}.
  \end{equation*}
\end{lemma}

Finally, we recall the Gagliardo--Nirenberg inequalities in Lorentz spaces. They will be used to obtain Sobolev inequalities and interpolation inequalities.

\begin{lemma}[{Gagliardo--Nirenberg inequality in Lorentz spaces \protect\cite[Theorem 2.1]{H-Y-Z}}] \label{lemma:GN}
  Let $1\leq p < \infty$, $1\leq p_2, q, q_1, q_2<\infty$, $0<\alpha<q$, and assume that $0<s<N$, $1<p_1<\frac{N}{s}$. If 
  \begin{gather*}
    \frac{\alpha}{q_1} + \frac{q-\alpha}{q_2} = 1, \quad \text{and} \quad
    \alpha \left( \frac{1}{p_1} - \frac{s}{N} \right) + (q - \alpha) \frac{1}{p_2} = \frac{q}{p},
  \end{gather*}
  then there exists $A>0$ such that
  \begin{equation*}
    \norm{u}_{L^{p,q}(\R^N)} \leq A \norm{D^{s}u}_{L^{p_1,q_1}(\R^N)}^{\frac{\alpha}{q}} \norm{u}_{L^{p_2,q_2}(\R^N)}^{\frac{q-\alpha}{q}},
  \end{equation*}
  where $\omega_{N-1}$ denotes the surface measure of the unit sphere $S^{N-1}\subset\R^N$ and
  \begin{equation*}
    A \leq \left( \frac{3Np_1}{N - p_1s} \left( \frac{\omega_{N-1}}{N} \right)^{\frac{N-s}{N}} \right)^{\frac{\alpha}{q}}.
  \end{equation*}
\end{lemma}

\begin{lemma}[{Gagliardo--Nirenberg inequality in Lorentz spaces \protect\cite[Theorem 5.10]{B-K-O}}] \label{lemma:GN2}
  Let $s, s_1, s_2\in\R$, $1<p,p_1,p_2<\infty$, $1\leq q,q_1,q_2\leq\infty$, and $0<\theta<1$. Assume that  
  \begin{equation*}
    (1-\theta)s_1 + \theta s_2 - s = N\left( \frac{1-\theta}{p_1} + \frac{\theta}{p_2} - \frac{1}{p} \right) \geq 0.
  \end{equation*}
  If at least one of the following conditions holds:
  \begin{enumerate}[label=(\roman*)]
    \item $\frac{1-\theta}{q_1} + \frac{\theta}{q_2} \geq \frac{1}{q}$, 
    \item $s=s_1=s_2$, $p_1\neq p_2$, 
    \item $(1-\theta)s_1 + \theta s_2 > s$, $s_2 - s_1 \neq \frac{N}{p_2} - \frac{N}{p_1}$, 
  \end{enumerate}
  then, for every $f$ satisfying $D^{s_1}f\in L^{p_1,q_1}$, $D^{s_2}f\in L^{p_2,q_2}$, the interpolation inequality 
  \begin{equation*}
    \norm{D^{s} f}_{L^{p,q}} \lesssim \norm{D^{s_1} f}_{L^{p_1,q_1}}^{1-\theta} \norm{D^{s_2} f}_{L^{p_2,q_2}}^{\theta}
  \end{equation*}
  holds. 
\end{lemma}

\begin{remark} \label{remark:2-8}
  \refsp{\lem:GN2} gives the Sobolev-type inequality 
  \begin{equation*}
    \norm{f}_{L^{p,q}} \lesssim \norm{D^{\sigma} f}_{L^{p_1,q}}, \quad 1<p_1<\frac{N}{\sigma}, \quad 1\leq q\leq\infty. 
  \end{equation*}
  Indeed, taking 
  \begin{equation*}
    p_2 = p, \quad q_1 = q_2 = q, \quad s_1 = \sigma, \quad s = s_2 = 0, \quad \frac{1}{p} = \frac{1}{p_1} - \frac{s}{N},
  \end{equation*}
  we see that the assumptions of \refsp{\lem:GN2} are satisfied.
\end{remark}

\subsection{Strichartz estimates} \label{sec:2-3}

Strichartz estimates originate from the Fourier restriction theorem of Strichartz~\cite{str}, and are widely used as space-time estimates for solutions to dispersive equations. For the Schr\"odinger equation, he first established the estimate in the case
\begin{equation*}
  q = r = \frac{2(N+2)}{N}.
\end{equation*}
The range of exponents was later extended by Ginibre--Velo~\cite{G-V}, Yajima~\cite{yaj}, and others. Keel--Tao~\cite{K-T} further established Strichartz estimates including the endpoint case in an abstract framework based on the energy and dispersive estimates. In this subsection, we summarize the Strichartz estimates in Lorentz spaces. We first recall the definition of admissible pairs in the form needed below.

\begin{definition} \label{definition:adm}
    Generalizing the notation of admissible pairs introduced in \refsp{sec:1-2}, for $s\in\R$, a pair $(q,r)$ satisfying
    \begin{equation*}
        \frac{2}{q} + \frac{N}{r} = \frac{N}{2} - s
    \end{equation*} 
    is called an $\dot{H}^s$-admissible pair. For $s<\frac{N}{2}$, we define
    \begin{equation*}
        \cA_s \coloneq \set{(q,r):\dot{H}^s\text{-admissible pair}}{\frac{2N}{N-2s}<r<\frac{2N}{N-2}}, 
    \end{equation*}
    where, for $N=1,2$, the condition on $r$ is replaced by $\frac{2N}{N-2s}<r<\infty$. In particular, when $s=0$, this definition agrees with $\cA_0$ introduced in \refsp{sec:1-2}. 
\end{definition}

\begin{proposition}[{\kern-3pt\cite[Theorem 10.1]{K-T}}] \label{proposition:2}
  \begin{enumerate}[label=(\roman*)]
    \item Let $(q,r)$ be an admissible pair with $r<\infty$. Then there exists a constant $C>0$ such that, for every $\varphi\in L^2$,
    \begin{equation} \label{eq:14}
      \norm{e^{it\Delta}\varphi}_{L^q(\R; L^{r,2})} \leq C\lnorm{\varphi}.
    \end{equation}
    \item Let $(q_i,r_i)$, where $i=1,2$, be two admissible pairs with $r_1, r_2<\infty$. Then there exists a constant $C>0$ such that, for every $f\in L^{q_2^\prime}(\R; L^{r_2^\prime, 2})$, we have
    \begin{equation} \label{eq:15}
      \norm{\int_{0}^{t} e^{i(t-s)\Delta}f(s,\cdot)ds}_{L^{q_1}(\R; L^{r_1,2})} \leq C\norm{f}_{L^{q_2^\prime}(\R; L^{r_2^\prime, 2})}.
    \end{equation}
  \end{enumerate}
\end{proposition}

The inhomogeneous Strichartz estimate stated below is known to hold for a wider range of exponents than the homogeneous one. Results in this direction go back to Cazenave--Weissler~\cite{C-W} and Kato~\cite{kat}, and a wider range of exponents for the Schr\"odinger equation was later obtained by Foschi~\cite{fos} and Vilela~\cite{vil}. Taggart~\cite{tag} further generalized these results to the abstract framework of Keel--Tao. In this paper, we use this extension of the Strichartz estimate in the proof of continuous dependence in the critical case. Here, we introduce the following definition, which corresponds to the case $\sigma = \frac{N}{2}$ in the definition of a $\sigma$-acceptable pair in Taggart~\cite{tag}.

\begin{definition} \label{definition:2-6}
  We call $(q,r)$ an acceptable pair if 
  \begin{equation*}
    \frac{1}{q} < \frac{N}{2} - \frac{N}{r}, \hspace{20pt} 1\leq q < \infty,\ 2<r<\infty
  \end{equation*}
  or if $(q,r) = (\infty,2)$. 
\end{definition}

\begin{proposition}[{\kern-3pt\cite[Theorem 1.4, p. 828]{tag}}] \label{proposition:2-12}
  Let $N\geq 1$. Suppose that $(q,r)$ and $(\widetilde{q}, \widetilde{r})$ are two acceptable pairs satisfying the scaling condition
  \begin{equation*}
    \frac{1}{q} + \frac{1}{\widetilde{q}} = \frac{N}{2}\left( 1-\frac{1}{r}-\frac{1}{\widetilde{r}} \right),
  \end{equation*}
  and
  \begin{equation*}
    \frac{1}{q} + \frac{1}{\widetilde{q}} < 1, \hspace{10pt} \frac{N-2}{r} \leq \frac{N}{\widetilde{r}}, \hspace{10pt} \frac{N-2}{\widetilde{r}} \leq \frac{N}{r}.
  \end{equation*}
  Then there exists a constant $C>0$ such that, for every $f\in L^{\widetilde{q}^\prime}(\R; L^{\widetilde{r}^\prime, 2})$, we have
    \begin{equation} \label{eq:16}
      \norm{\int_{0}^{t} e^{i(t-s)\Delta}f(s,\cdot)ds}_{L^{q}(\R; L^{r,2})} \leq C\norm{f}_{L^{\widetilde{q}^\prime}(\R; L^{\widetilde{r}^\prime, 2})}.
    \end{equation}
\end{proposition}

\begin{remark} \label{remark:4}
  The inhomogeneous Strichartz estimates \refsp{eq:15} and \refsp{eq:16} remain valid when the time interval $\R$ is replaced by any bounded or unbounded interval in $\R$. 
\end{remark}

\subsection{Classes of admissible coefficients} \label{sec:2-4}

In this subsection, we specify the function spaces to which the nonlinear coefficient $a(x)$ belongs. For two function spaces $X$ and $Y$, we write 
\begin{equation*}
    X + Y \coloneq \set{f = f_1 + f_2}{f_1\in X, f_2\in Y}.
\end{equation*}
When necessary, we equip this space with the sum-space norm
\begin{equation*}
    \norm{f}_{X+Y} = \inf_{f = f_1 + f_2} \left( \norm{f_1}_{X} + \norm{f_2}_{Y} \right).
\end{equation*}
Throughout this subsection, we fix $0\leq s\leq 1$. 

\begin{assumption*} \phantomsection \label{assumption:coeff}
    Assume that there exists $0<b<N-s$ such that
    \begin{equation} \label{eq:2-7}
        a \in L^{\infty} + L^{\frac{N}{b},\infty},
    \end{equation}
    where $L^{\frac{N}{b},\infty} = L^{\infty}$ when $b=0$. When $s>0$, we further assume that
    \begin{equation} \label{eq:2-8}
        D^s a \in L^{\infty} + L^{\frac{N}{b+s},\infty}.
    \end{equation}
    When $s=0$, we assume only \refsp{eq:2-7}. 
\end{assumption*}

The exponent $b$ is fixed in the following arguments.

\begin{remark} \label{remark:2-15}
    For $0\leq \tau\leq s$ and $b_{\ast}\in \{0,b\}$, we define 
    \begin{equation*}
        p_{\tau, b_{\ast}} \coloneq \begin{dcases}
            \frac{N}{b_{\ast} + \tau}, & b_{\ast} + \tau > 0, \\
            \infty, & b_{\ast} + \tau = 0.
        \end{dcases}
    \end{equation*}
    This notation will be used frequently in the subsequent nonlinear estimates.
    
    With this notation, under the \assumptionref, we decompose $a$ and $D^s a$ separately as
    \begin{equation*}
        a = a_{\infty} + a_b, \qquad a_{\infty}\in L^{\infty} , \quad a_{b} \in L^{p_{0,b},\infty},
    \end{equation*} 
    and
    \begin{equation*}
        D^s a = A_{\infty} + A_{b}, \qquad A_{\infty}\in L^{\infty} , \quad A_{b} \in L^{p_{s,b},\infty}.
    \end{equation*} 
    Here, the two decompositions are independent, and we do not require $D^sa_b = A_b$. 

    In the nonlinear estimates, we apply the fractional Leibniz rule to the full coefficient $a$, and then estimate the term involving $a_{\infty}$, $a_{b}$, $A_{\infty}$, and $A_{b}$ separately. The terms with $A_{\infty}$ are estimated using the boundedness of this component. The terms involving the $L^\infty$-components $a_{\infty}$ and $A_{\infty}$ of $a$ and its derivatives are treated separately by standard estimates as in the classical NLS case. We omit the details of the estimates for $a_{\infty}$ and $A_{\infty}$.
\end{remark}

\begin{remark}\label{remark:2-15.5}
    Under the \assumptionref, for every $0<\sigma<s$,
    \begin{equation*}
        D^{\sigma}a\in L^{\infty} + L^{\frac{N}{b+\sigma},\infty}.
    \end{equation*}

    Indeed, let $\chi\in C_0^\infty(\R^N)$ be a smooth cutoff function with $0\leq\chi\leq 1$, $\chi(\xi) = 1$ for $|\xi|\leq 1$. Then we can write
    \begin{align*}
        D^{\sigma}a &= \F^{-1}\left[|\xi|^{\sigma}\chi(\xi)\widehat{a}\right] + \F^{-1}\left[|\xi|^{\sigma-s}(1-\chi(\xi)) |\xi|^{s} \widehat{a}\right] \\
        &= \chi_1\ast a + \chi_2\ast D^s a,
    \end{align*}
    where $\chi_1 = \F^{-1}[|\xi|^{\sigma}\chi(\xi)]$ and $\chi_2 = \F^{-1}[|\xi|^{\sigma-s}(1-\chi(\xi))]$. By a dyadic decomposition, we obtain
    \begin{equation*}
        |\chi_1(x)|\lesssim \jb{x}^{-N-\sigma},
    \end{equation*}
    and 
    \begin{equation*}
        |\chi_2(x)| \lesssim \begin{dcases}
            |x|^{-N+s-\sigma}, & 0<|x|\leq 1, \\
            |x|^{-N-s}, & |x|> 1.
        \end{dcases}
    \end{equation*}
    Since $0<\sigma<s$, these estimates imply that
    \begin{equation*}
        \chi_1\in L^1, \quad \chi_2\in L^1(\R^N) \cap L^{\frac{N}{N-s+\sigma},\infty}.
    \end{equation*}
    Thus, by Young's inequality, we obtain
    \begin{align*}
        \norm{\chi_1\ast a}_{L^{\infty} + L^{\frac{N}{b+\sigma},\infty}} &\lesssim \norm{\chi_1\ast a}_{L^{\infty} + L^{\frac{N}{b},\infty}} \\
        &\lesssim \norm{\chi_1}_{L^1}\norm{a}_{L^{\infty} + L^{\frac{N}{b},\infty}} \lesssim \norm{a}_{L^{\infty} + L^{\frac{N}{b},\infty}}.
    \end{align*}

    Next, we use O'Neil's inequality (see \cite[T\textsc{heorem} 2.6]{one}), which states that 
    \begin{equation*}
        \norm{f\ast g}_{L^{r,\infty}} \lesssim \norm{f}_{L^{p,\infty}} \norm{g}_{L^{q,\infty}},
    \end{equation*}
    where $1<p,q,r<\infty$ and $\frac{1}{p} + \frac{1}{q} = 1 + \frac{1}{r}$. Since $1 + \frac{b+\sigma}{N} = \frac{N-s+\sigma}{N} + \frac{b+s}{N}$, by \refsp{eq:2-8} and O'Neil's inequality, we obtain
    \begin{align*}
        \norm{\chi_2\ast D^s a}_{L^{\infty} + L^{\frac{N}{b+\sigma},\infty}} \lesssim \norm{D^s a}_{L^{\infty} + L^{\frac{N}{b+s},\infty}}.
    \end{align*}
\end{remark}

\begin{remark}\label{remark:2-16}
    First, let $a(x) = \mu$. Then $a\in L^{\infty}$ and $D^s a = 0$ when $s>0$. Thus, $a$ satisfies the \assumptionref\ with $b=0$. In this case, equation \refsp{eq:1-1} reduces to the classical NLS with a power-type nonlinearity. 
    
    Next, let $a(x) = \lambda |x|^{-b}$ with $0<b<N-s$. Since
    \begin{equation*}
        d_{|x|^{-b}}(\tau) = \left| \set{x\in\R^N}{|x|^{-b} > \tau} \right| = C_N \tau^{-\frac{N}{b}},
    \end{equation*}
    we have $a\in L^{\frac{N}{b},\infty}$. Moreover, 
    \begin{equation*}
        D^s |x|^{-b} = \F^{-1}\left[ |\xi|^{s} \widehat{|x|^{-b}} \right] = C_{N,b} \F^{-1}\left[ |\xi|^{s} |\xi|^{b-N} \right] = C_{N,b,s} |x|^{-b-s}.
    \end{equation*}
    Therefore, $D^s a \in L^{\frac{N}{b+s},\infty}$ and $a$ is covered by \assumptionref.

    The coefficient $a(x) = \lambda |x|^{-b} + \mu$ falls within the class given in \assumptionref. Indeed, 
    \begin{equation*}
        a\in L^{\frac{N}{b},\infty} + L^{\infty}, \quad D^s a \in L^{\frac{N}{b+s}, \infty} \quad (s>0).
    \end{equation*}

    Furthermore, consider
    \begin{equation*}
        a(x) = \mu + \sum_{j=1}^{\infty} \lambda_j |x - y_j|^{-b_j},
    \end{equation*}
    where $0<b_j\leq b<N-s$ and $\{\lambda_j\}_{j=1}^{\infty} \in \ell^1$. For $0\leq b_j\leq b$, we have 
    \begin{gather*}
        \norm{|x - y_j|^{-b_j}}_{L^{\infty} + L^{\frac{N}{b},\infty}} \lesssim 1, \\
        \norm{D^s |x - y_j|^{-b_j}}_{L^{\infty} + L^{\frac{N}{b+s},\infty}} \lesssim 1.
    \end{gather*}
    where the constants are independent of $j$, $b_j$, and $y_j$. Thus, since $\{\lambda_j\} \in \ell^1$, the series converges in the sum spaces. In particular,
    \begin{equation*}
        a \in L^{\infty} + L^{\frac{N}{b},\infty}, \quad D^s a \in L^{\infty} + L^{\frac{N}{b+s},\infty}.
    \end{equation*}
    Hence, $a$ satisfies the \assumptionref.
    \end{remark}

\section{ Difference estimates in Lorentz spaces} \label{sec:3}

In this section, we prepare difference estimates for the nonlinear term that will be used to prove continuous dependence on the initial data in the main theorems. In particular, when $0<\alpha<1$, the derivative of the nonlinear function is only H\"older continuous. Therefore, the usual Lipschitz-type difference estimate is not sufficient, and a more delicate analysis is required.

Fujiwara~\cite{fuj} established fractional derivative bounds for the difference of such nonlinear mappings by using the Littlewood--Paley decomposition and pointwise difference estimates. In this section, we extend the argument of Fujiwara~\cite{fuj} to the setting of Lorentz spaces and establish the nonlinear bounds needed in the proofs of the main theorems. In \refsp{sec:3-1}, we prepare some basic facts about the Littlewood--Paley decomposition and the Hardy--Littlewood maximal operator. In \refsp{sec:3-2}, we prove the difference estimates that forms the main result of this section.

\subsection{Preliminary estimates} \label{sec:3-1}

In this subsection, we collect some basic facts on the Littlewood--Paley decomposition and the Hardy--Littlewood maximal operator that will be used in the proof of the difference estimates.

First, we choose a radial function $\psi\in\Sc(\R^N)$ such that
\begin{equation*}
  \supp{\widehat{\psi}} \subset \set{\xi\in\R^N}{\frac{1}{2} < |\xi| < 2}, \quad 0 \leq \widehat{\psi} \leq 1,
\end{equation*}
and
\begin{equation*}
  \sum_{j\in\Z} \widehat{\psi}(2^{-j}\xi) = 1 \quad (\xi\neq 0).
\end{equation*}
We further define
\begin{equation*}
  \psi_j(x) \coloneq 2^{jN}\psi(2^jx), \quad Q_jf \coloneq \psi_j \ast f,
\end{equation*}
and
\begin{equation*}
  \widetilde{Q}_jf \coloneq Q_{j-1}f + Q_jf + Q_{j+1}f.
\end{equation*}
Then $\widetilde{Q}_jQ_j = Q_j$. 

Throughout this section, for $0<r<\infty$, we use the notation 
\begin{equation} \label{eq:3-A}
  \norm{\{f_j\}}_{L^{p,q}(\ell^r)} \coloneq \norm{\left(\sum_{j\in\Z} |f_j|^r\right)^{\frac{1}{r}}}_{L^{p,q}}. 
\end{equation}
For $r=\infty$, we use the usual modification. Then, for $s\in\R$, $1<p<\infty$, $1\leq q\leq \infty$, and suitable $f$,
\begin{equation} \label{eq:3-B}
  \norm{D^s f}_{L^{p,q}} \sim \norm{\{2^{js}Q_jf\}_{j\in\Z}}_{L^{p,q}(\ell^2)}.
\end{equation}
See \cite[Theorem 5.1]{B-K-O} and the references therein. 

The Hardy--Littlewood maximal operator $M$ is defined by
\begin{equation*}
  Mf(x) \coloneq \sup_{r>0} \frac{1}{|B_r(x)|}\int_{B_r(x)} |f(y)| dy.
\end{equation*}
For $a>0$, define
\begin{equation*}
  M^{(a)}f(x) \coloneq \left( M\left(|f|^a\right)(x) \right)^{\frac{1}{a}}.
\end{equation*}
Since $M:L^1\to L^{1,\infty}$ is bounded (see \cite[Theorem 2.1.6]{gra}), interpolation yields
\begin{equation*}
  \norm{Mf}_{L^{p,q}} \lesssim \norm{f}_{L^{p,q}}, \quad 1<p<\infty, \quad 0<q\leq \infty
\end{equation*}
Thus, $M$ is bounded on $L^{p,q}$. Moreover, by \cite[Lemma 4.5 (i)]{L-Y-Y}, the vector-valued estimate 
\begin{equation*}
  \norm{\{Mf_j\}}_{L^{p,q}(\ell^r)} \lesssim \norm{\{f_j\}}_{L^{p,q}(\ell^r)}
\end{equation*}
also holds for $1<r\leq\infty$. 

We begin by recalling two lemmas that will be used together with the Littlewood--Paley decomposition.

\begin{lemma}[{\kern-4.5pt\cite[Theorem 2.1.10]{gra}}] \label{lemma:1}
  Let $w\in L^1$ be a positive, radial, and decreasing function. Then, for every $x\in\R^N$,
  \begin{equation*}
    |w\ast g(x)| \leq \lpnorm{1}{w} Mg(x).
  \end{equation*}
\end{lemma}

\begin{lemma}[{\kern-3pt\cite[Lemma 2]{fuj}}] \label{lemma:2}
  For every $k\in\Z$ and $x\in\R^N$, 
  \begin{equation*}
    \left|\widetilde{Q}_kg(x)\right| \lesssim Mg(x).
  \end{equation*}
  Moreover, if $|y| < 2^{-k}$, then for every $x\in\R^N$,
  \begin{equation*}
    \left| \widetilde{Q}_k g(x+y) - \widetilde{Q}_k g(x) \right| \lesssim 2^{k}|y|Mg(x),
  \end{equation*}
  and
  \begin{equation*}
    \left| \widetilde{Q}_k g(x+y) - 2\widetilde{Q}_k g(x) + \widetilde{Q}_k g(x-y) \right| \lesssim 2^{2k}|y|^2Mg(x).
  \end{equation*}
\end{lemma}

We also need the following estimate to handle the dyadic sums arising from the Littlewood--Paley decomposition.

For $c\in\R$, we define $c_{-} \coloneq \min\{c,0\}$. 

\begin{lemma}[{\kern-3pt\cite[Lemma 3]{fuj}}] \label{lemma:3}
  Let $0<s<s^{\prime}$, $q\geq 1$, and $a\in\ell^q$. Then
  \begin{equation*}
    \norm{2^{js}\sum_{k}2^{(k-j)_{-}s^{\prime}}a_k}_{\ell_j^q} \leq \left( \frac{2^s}{2^s - 1} + \frac{1}{2^{s^{\prime} - s} - 1} \right) \norm{2^{ks}a_k}_{\ell_k^q}.
  \end{equation*}
\end{lemma}

The previous lemmas provide the bounds that will be used below. We next recall a basic estimate of Fujiwara that plays an important role in the difference argument.

\begin{lemma}[{\kern-3pt\cite[Lemma 4]{fuj}}] \label{lemma:4}
  Let $0<s<1$ and $1\leq p, q < \infty$. Then
  \begin{align*}
    &\norm{2^{js}\norm{v(y+x)\left( u(y+x) - u(x) \right)|\psi_j|^{\frac{1}{p}}}_{L_y^p}}_{\ell_j^q} \\
    &\hspace{80pt} \lesssim M^{(p)} v(x) \norm{2^{ks}MQ_ku}_{\ell_k^q} + \norm{2^{ks}M^{(p)}(vMQ_ku)(x)}_{\ell_k^q}.
  \end{align*}
\end{lemma}

Using \refsp{\lem:3} and the maximal function estimates in Lorentz spaces, we obtain the following Lorentz-space version of \cite[Corollary 2]{fuj}.

\begin{corollary} \label{corollary:2}
  Let $0<s<1$, $1\leq p_0 < 2$, and $p_0 < p < \infty$. Suppose that $1<p_1,p_2<\infty$ and $1\leq q, q_1, q_2\leq \infty$ satisfy
  \begin{equation*}
    \frac{1}{p} = \frac{1}{p_1} + \frac{1}{p_2}, \hspace{10pt} \frac{1}{q} \leq \frac{1}{q_1} + \frac{1}{q_2}.
  \end{equation*}
  Then
  \begin{equation*}
    \norm{\left\{2^{js}\norm{\left( u(y+x) - u(x) \right)|\psi_j|^{\frac{1}{p_0}}}_{L_y^{p_0}}\right\}}_{L^{p,q}(\ell_j^2)} \lesssim \norm{D^su}_{L^{p,q}},
  \end{equation*}
  and
  \begin{align*}
    & \norm{\left\{2^{js}\norm{v(y+x) \left( u(y+x) - u(x) \right)|\psi_j|^{\frac{1}{p_0}}}_{L_y^{p_0}}\right\}}_{L^{p,q}(\ell_j^2)} \\
    & \hspace{180pt} \lesssim \norm{v}_{L^{p_1,q_1}} \norm{D^su}_{L^{p_2,q_2}}.
  \end{align*}
\end{corollary}

\begin{proof}
  First, the assumption $p>p_0$ and the relation between $p, p_1, p_2$ imply that $p_0 < p_1$, $p_0 < p_2$. We also note that 
  \begin{equation*}
    \norm{M^{(p_0)}f}_{L^{p_1,q_1}} = \norm{M|f|^{p_0}}_{L^{\frac{p_1}{p_0}, \frac{q_1}{p_0}}}^{\frac{1}{p_0}} \lesssim \norm{|f|^{p_0}}_{L^{\frac{p_1}{p_0}, \frac{q_1}{p_0}}}^{\frac{1}{p_0}} = \norm{f}_{L^{p_1,q_1}}.
  \end{equation*}
  Applying \refsp{\lem:4} with $v\equiv 1$, $p=p_0$, and $q=2$, taking the $L^{p,q}$-norm, and using the vector-valued maximal function estimate together with \refsp{eq:3-A} and \refsp{eq:3-B}, we obtain
  \begin{align*}
    &\norm{\left\{2^{js}\norm{\left( u(y+x) - u(x) \right)|\psi_j|^{\frac{1}{p_0}}}_{L_y^{p_0}}\right\}}_{L^{p,q}(\ell_j^2)} \\
    &\hspace{30pt} \lesssim \norm{\left\{2^{ks}MQ_ku\right\}}_{L^{p,q}(\ell_k^2)} + \norm{\left\{2^{ks}M^{(p_0)}(MQ_ku)(x)\right\}}_{L^{p,q}(\ell_k^2)} \\
    &\hspace{30pt} \lesssim \norm{\left\{2^{ks}Q_ku\right\}}_{L^{p,q}(\ell_k^2)} \sim \norm{D^su}_{L^{p,q}}.
  \end{align*}
  This proves the first estimate.

  Next, we apply \refsp{\lem:4} with $p=p_0$ and $q=2$. Taking the $L^{p,q}$-norm yields
  \begin{align*}
    & \norm{\left\{2^{js}\norm{v(y+x) \left( u(y+x) - u(x) \right)|\psi_j|^{\frac{1}{p_0}}}_{L_y^{p_0}}\right\}}_{L^{p,q}(\ell_j^2)} \\
    & \hspace{60pt} \lesssim \norm{M^{(p_0)} v(x) \norm{2^{ks}MQ_ku}_{\ell_k^2}}_{L^{p,q}} \\
    & \hspace{120pt} + \norm{\left\{2^{ks}M^{(p_0)}(vMQ_ku)(x)\right\}}_{L^{p,q}(\ell_k^2)} \\
    & \hspace{60pt} \lesssim \norm{M^{(p_0)}v}_{L^{p_1,q_1}} \norm{\left\{2^{ks}MQ_ku\right\}}_{L^{p_2,q_2}(\ell_k^2)} \\
    & \hspace{120pt} + \norm{|v|\norm{2^{ks}MQ_ku}_{\ell_k^2}}_{L^{p,q}} \\
    & \hspace{60pt} \lesssim \norm{v}_{L^{p_1,q_1}} \norm{\left\{2^{ks}Q_ku\right\}}_{L^{p_2,q_2}(\ell_k^2)} \\
    & \hspace{120pt} + \norm{v}_{L^{p_1,q_1}} \norm{\left\{2^{ks}MQ_ku\right\}}_{L^{p_2,q_2}(\ell_k^2)} \\
    & \hspace{60pt} \lesssim \norm{v}_{L^{p_1,q_1}} \norm{D^su}_{L^{p_2,q_2}}.
  \end{align*}
  This proves the second estimate.
\end{proof}

\subsection{ Nonlinear difference estimates} \label{sec:3-2}

In this subsection, we derive estimates for $D^s(F(u)-F(v))$ using the results established in the previous subsection. The following proposition extends the estimates of Fujiwara~\cite{fuj} to the setting of Lorentz spaces, and will be used later in the proof of continuous dependence on the initial data. 

\begin{proposition} \label{proposition:Fujiwara}
  Let $N\geq 1$, $\alpha>0$, and $0<s<1$. Suppose that $F\in C^1(\C)$ satisfies $F(0) = F^\prime(0) = 0$ and, for any $u,v\in\C$,
  \begin{equation*}
    \left|F(u) - F(v)\right| \lesssim (|u|+|v|)^{\alpha} |u-v|, 
  \end{equation*}
  and
  \begin{equation}
    \left|F^\prime(u) - F^\prime(v)\right| \lesssim \begin{dcases}
      (|u|+|v|)^{\alpha-1} |u-v|, & \alpha\geq 1, \\
      |u-v|^{\alpha}, & \alpha<1.
    \end{dcases} \label{eq:3-2}
  \end{equation}
  Here, $F^{\prime}$ denotes the real derivative of $F$ when $\C$ is identified with $\R^2$. Let $1<p,p_1,p_2<\infty$ and $1\leq q,q_1,q_2\leq \infty$ satisfy
  \begin{equation*}
    \frac{1}{p} = \frac{1}{p_1} + \frac{1}{p_2}, \quad \frac{1}{q} \leq \frac{1}{q_1} + \frac{1}{q_2}
  \end{equation*}
  and assume that $\alpha p_1>1$. Suppose further that 
  \begin{equation*}
    u,v\in L^{\alpha p_1, \alpha q_1}, \quad \text{and} \quad D^s u, D^s v\in L^{p_2, q_2}.
  \end{equation*}
  If $\alpha\geq 1$, then
  \begin{align}
    &\norm{D^s(F(u) - F(v))}_{L^{p,q}} \notag\\
    &\hspace{20pt} \lesssim \left( \norm{u}_{L^{\alpha p_1, \alpha q_1}} + \norm{v}_{L^{\alpha p_1, \alpha q_1}} \right)^{\alpha} \norm{D^s(u-v)}_{L^{p_2, q_2}} \notag\\
    & \hspace{40pt} + \left( \norm{D^su}_{L^{p_2, q_2}} + \norm{D^sv}_{L^{p_2, q_2}} \right) \notag\\
    & \hspace{60pt} \times \left( \norm{u}_{L^{\alpha p_1, \alpha q_1}} + \norm{v}_{L^{\alpha p_1, \alpha q_1}} \right)^{\alpha-1} \norm{u-v}_{L^{\alpha p_1, \alpha q_1}}. \label{eq:3-3}
  \end{align}

  On the other hand, if $0<\alpha<1$, then for any $0<\theta\leq \alpha$, 
  \begin{align}
    &\norm{D^s(F(u) - F(v))}_{L^{p,q}} \notag\\
    &\hspace{20pt} \lesssim \left( \norm{u}_{L^{\alpha p_1, \alpha q_1}} + \norm{v}_{L^{\alpha p_1, \alpha q_1}} \right)^{\alpha} \norm{D^s(u-v)}_{L^{p_2, q_2}} \notag\\
    &\hspace{40pt} + \left( \norm{D^su}_{L^{p_2, q_2}} + \norm{D^sv}_{L^{p_2, q_2}} \right) \notag\\
    &\hspace{60pt} \times \left( \norm{u}_{L^{\alpha p_1, \alpha q_1}} + \norm{v}_{L^{\alpha p_1, \alpha q_1}} \right)^{\alpha-\theta} \norm{u-v}_{L^{\alpha p_1, \alpha q_1}}^{\theta}. \label{eq:3-4}
  \end{align}
\end{proposition}

\begin{proof}
  Since $\supp{\widehat{\psi}}$ does not contain the origin, 
  \begin{equation*}
    \int_{\R^N} \psi_j(y) dy = 0.
  \end{equation*}
  Hence, 
  \begin{align*}
    & Q_j (F(u) - F(v))(x) \\
    & \hspace{20pt} = \int_{\R^N} \left[ F(u(x+y)) - F(u(x)) - F(v(x+y)) + F(v(x)) \right] \psi_j(y) dy.
  \end{align*}
  By the Littlewood--Paley characterization \refsp{eq:3-B},
  \begin{align*}
    &\norm{D^s(F(u) - F(v))}_{L^{p,q}}\\
    &\hspace{20pt} \lesssim \Bigg\|\Bigg\{2^{js}\int_{\R^N} \left| F(u(x+y)) - F(u(x)) - F(v(x+y)) + F(v(x)) \right|\\
    &\hspace{200pt} \times |\psi_j(y)| dy\Bigg\}\Bigg\|_{L^{p,q}(\ell^2)}. 
  \end{align*}
  By the fundamental theorem of calculus and \refsp{eq:3-2}, we have
  \begin{align}
    & |F(u(x+y)) - F(u(x)) - F(v(x+y)) + F(v(x))| \notag \\
    & \hspace{20pt} \lesssim \left( |u(x+y)|^{\alpha} + |u(x)|^{\alpha} \right)|(u-v)(x+y) - (u-v)(x)| \notag \\
    & \hspace{40pt} + 
    \begin{dcases}
      \left( |u(x+y)| + |u(x)| + |v(x+y)| + |v(x)| \right)^{\alpha - 1}  \\
      \hspace{20pt} \times \left( |(u-v)(x+y)| + |(u-v)(x)| \right) \\
      \hspace{40pt} \times |v(x+y) - v(x)|, \qquad \alpha\geq 1, \\
      \left( |(u-v)(x+y)|^{\alpha} + |(u-v)(x)|^{\alpha} \right) \\
      \hspace{40pt} \times |v(x+y) - v(x)|, \qquad 0<\alpha<1.
    \end{dcases} \label{eq:3-6}
  \end{align}

  First, we consider the case $0<\alpha<1$. For the first term in \refsp{eq:3-6}, applying \refsp{\cor:2} with $p_0=1$, we obtain
  \begin{align*}
    & \norm{\left\{2^{js}\int_{\R^N} \left| |u(\ \cdot+y)|^{\alpha} |(u-v)(\ \cdot+y) - (u-v)(\ \cdot\ )| \right| |\psi_j(y)| dy\right\}}_{L^{p,q}(\ell^2)} \\
    & \hspace{20pt} \lesssim \norm{|u|^{\alpha}}_{L^{p_1,q_1}} \norm{D^s(u-v)}_{L^{p_2,q_2}} \\
    & \hspace{20pt} \lesssim \norm{u}_{L^{\alpha p_1,\alpha q_1}}^{\alpha} \norm{D^s(u-v)}_{L^{p_2,q_2}}.
  \end{align*}
  A similar estimate holds for the term with $|u(x)|^{\alpha}$ in place of $|u(x+y)|^{\alpha}$. 

  For the second term in \refsp{eq:3-6}, another application of \refsp{\cor:2} gives
  \begin{align*}
    & \Bigg\|\Bigg\{2^{js}\int_{\R^N} \left( |(u-v)(\ \cdot+y)|^{\alpha} + |(u-v)(\ \cdot\ )|^{\alpha} \right) \\
    & \hspace{100pt} \times |v(\ \cdot +y) - v(\ \cdot\ )| |\psi_j(y)| dy\Bigg\}\Bigg\|_{L^{p,q}(\ell^2)} \\
    & \hspace{20pt} \lesssim \norm{u-v}_{L^{\alpha p_1,\alpha q_1}}^{\alpha} \norm{D^sv}_{L^{p_2,q_2}}.
  \end{align*}
  Symmetry between $u$ and $v$ then yields
  \begin{align*}
    &\norm{D^s(F(u) - F(v))}_{L^{p,q}} \notag\\
    &\hspace{20pt} \lesssim \left( \norm{u}_{L^{\alpha p_1, \alpha q_1}} + \norm{v}_{L^{\alpha p_1, \alpha q_1}} \right)^{\alpha} \norm{D^s(u-v)}_{L^{p_2, q_2}} \notag\\
    &\hspace{40pt} + \left( \norm{D^su}_{L^{p_2, q_2}} + \norm{D^sv}_{L^{p_2, q_2}} \right) \norm{u-v}_{L^{\alpha p_1, \alpha q_1}}^{\alpha}.
  \end{align*}
  Finally, by the quasi-triangle inequality in Lorentz spaces, 
  \begin{align*}
    \norm{u-v}_{L^{\alpha p_1,\alpha q_1}}^{\alpha} & = \norm{u-v}_{L^{\alpha p_1,\alpha q_1}}^{\alpha - \theta} \norm{u-v}_{L^{\alpha p_1,\alpha q_1}}^{\theta} \\
    & \lesssim \left( \norm{u}_{L^{\alpha p_1, \alpha q_1}} + \norm{v}_{L^{\alpha p_1, \alpha q_1}} \right)^{\alpha - \theta} \norm{u-v}_{L^{\alpha p_1,\alpha q_1}}^{\theta}.
  \end{align*}
  Combining these estimates gives \refsp{eq:3-4}. 

  Next, we deal with the case $\alpha\geq 1$. The first term in \refsp{eq:3-6} can be estimated in the same way as above. For the second term, we split it into finitely many terms and estimate each of them separately. For example, when $\alpha>1$, one of the resulting terms is
  \begin{equation*}
    \left( |u(x)| + |v(x)| \right)^{\alpha - 1} |(u-v)(x+y)| |v(x+y) - v(x)|.
  \end{equation*}
  Set 
  \begin{equation*}
    \frac{1}{\widetilde{p}} \coloneq \frac{1}{\alpha p_1} + \frac{1}{p_2}, \quad \frac{1}{\widetilde{q}} \coloneq \min \left\{ \frac{1}{\alpha q_1} + \frac{1}{q_2}, 1 \right\}.
  \end{equation*}
  Then $1<\widetilde{p}\leq \infty$ and $1\leq\widetilde{q}\leq \infty$. By H\"older's inequality in Lorentz spaces and \refsp{\cor:2}, 
  \begin{align*}
    & \Bigg\|\Bigg\{2^{js}\int_{\R^N} \left( |u(\ \cdot\ )| + |v(\ \cdot\ )| \right)^{\alpha-1} |(u-v)(\ \cdot +y)| \\
    & \hspace{100pt} \times |v(\ \cdot +y) - v(\ \cdot\ )| |\psi_j(y)| dy\Bigg\}\Bigg\|_{L^{p,q}(\ell^2)} \notag \\
    & \hspace{20pt} \lesssim \norm{\left( |u| + |v| \right)^{\alpha-1}}_{L^{\frac{\alpha p_1}{\alpha - 1}, \frac{\alpha q_1}{\alpha - 1}}} \notag \\
    & \hspace{40pt} \times \norm{\left\{2^{js}\int_{\R^N} |(u-v)(\ \cdot +y)| |v(\ \cdot +y) - v(\ \cdot\ )| |\psi_j(y)| dy\right\}}_{L^{\widetilde{p},\widetilde{q}}(\ell^2)} \notag \\
    & \hspace{20pt} \lesssim \norm{|u| + |v|}_{L^{\alpha p_1,\alpha q_1}}^{\alpha-1} \norm{u-v}_{L^{\alpha p_1,\alpha q_1}} \norm{D^sv}_{L^{p_2,q_2}} \notag \\
    & \hspace{20pt} \lesssim \left( \norm{u}_{L^{\alpha p_1,\alpha q_1}} + \norm{v}_{L^{\alpha p_1,\alpha q_1}} \right)^{\alpha-1} \norm{u-v}_{L^{\alpha p_1,\alpha q_1}} \norm{D^sv}_{L^{p_2,q_2}},
  \end{align*}
  where 
  \begin{equation*}
    \frac{1}{p} = \frac{\alpha - 1}{\alpha p_1} + \frac{1}{\widetilde{p}}, \quad \frac{1}{q} = \frac{\alpha - 1}{\alpha q_1} + \frac{1}{\widetilde{q}}.
  \end{equation*}

  Similarly, by H\"older's inequality in Lorentz spaces and the second estimate in \refsp{\cor:2}, we have
  \begin{align*}
    & \Bigg\|\Bigg\{2^{js}\int_{\R^N} \left( |u(\ \cdot +y)| + |v(\ \cdot+y)| \right)^{\alpha-1} |(u-v)(\ \cdot +y)| \\
    & \hspace{100pt} \times |v(\ \cdot +y) - v(\ \cdot\ )| |\psi_j(y)| dy\Bigg\}\Bigg\|_{L^{p,q}(\ell^2)} \\
    & \hspace{20pt} \lesssim \norm{\left( |u| + |v| \right)^{\alpha-1}|u-v|}_{L^{p_1, q_1}} \norm{D^sv}_{L^{p_2,q_2}} \\
    & \hspace{20pt} \lesssim \left( \norm{u}_{L^{\alpha p_1,\alpha q_1}} + \norm{v}_{L^{\alpha p_1,\alpha q_1}} \right)^{\alpha-1} \norm{u-v}_{L^{\alpha p_1,\alpha q_1}} \norm{D^sv}_{L^{p_2,q_2}}.
  \end{align*}

  The remaining terms are handled by a similar argument using the two estimates in \refsp{\cor:2} as appropriate. The same bounds also hold after interchanging $u$ and $v$. Therefore, \refsp{eq:3-3} follows. 
\end{proof}

\section{Local well-posedness in the subcritical case} \label{sec:4}

In this section, we prove \refsp{\thm:1-1} in the subcritical case. Throughout this section, we assume the hypotheses of \refsp{\thm:1-1} and consider the case $0<s\leq1$. The case $s=0$ can be treated in the same way by omitting the estimates involving derivatives.

As mentioned in \refsp{sec:1-2}, the subcritical condition gives $s_c < s$. Choose $\eta>0$ sufficiently small and set $\widetilde{s}\coloneq (1-\eta)s$. We fix $\eta$ small enough so that $s_c < \widetilde{s} < s$ and all the exponents appearing below are well-defined.

In the following two lemmas, we consider the estimates corresponding to $b_{\ast}\in\{0,b\}$. Since $b_{\ast}\leq b$, we have $s_c(b_{\ast})\leq s_c$.

We consider the following condition on:
\begin{equation} \label{eq:4-1}
  \begin{aligned}
    \max\biggl\{ \frac{N-2}{2N}, \frac{s}{N}, \frac{N - 2b_{\ast} + 2\alpha\widetilde{s}}{2N(\alpha + 1)} \biggr\} < \frac{1}{\rho} < \min\biggl\{ &\frac{1}{2}, \frac{N - 2b_{\ast} + 2\alpha\widetilde{s}+2}{2N(\alpha + 1)}, \\
      & \hspace{20pt} \frac{\widetilde{s}}{N} + \frac{1}{\alpha}, \frac{N-b_{\ast}+\alpha\widetilde{s}}{N(\alpha + 1)}.
    \biggr\}.
  \end{aligned}
\end{equation}

\begin{lemma} \label{lemma:4-1}
  The following statements hold. 
  \begin{enumerate}[label=(\roman*)]
    \item There exists an admissible pair $(\gamma,\rho)$ satisfying \refsp{eq:4-1}. 
    \item For any such pair $(\gamma,\rho)$, define
    \begin{equation*}
      \frac{1}{r} = 1 - \frac{b_{\ast}}{N} - \frac{\alpha + 1}{\rho} + \frac{\alpha\widetilde{s}}{N}.
    \end{equation*} 
    Then there exists $q$ such that $(q,r)\in\cA_0$.
    \item For the pairs $(\gamma,\rho)$ and $(q,r)$ above, define
    \begin{equation} \label{eq:4-2}
      \frac{1}{r_2}\coloneq \frac{1}{\rho} - \frac{\widetilde{s}}{N}, \quad \frac{1}{r_3}\coloneq \frac{1}{\rho} - \frac{s}{N}, \quad \frac{1}{r_{0,b_{\ast}}}\coloneq \frac{\alpha}{r_2} + \frac{1}{\rho}, \quad \frac{1}{r_{s,b_{\ast}}}\coloneq \frac{\alpha}{r_2} + \frac{1}{r_3}.
    \end{equation}
    Then these exponents are well-defined and satisfy 
    \begin{equation*}
      \frac{1}{r^{\prime}} = \frac{b_{\ast}}{N} + \frac{1}{r_{0,b_{\ast}}} = \frac{b_{\ast}+s}{N} + \frac{1}{r_{s,b_{\ast}}}, 
    \end{equation*}
    and
    \begin{equation*}
      \frac{1}{q^{\prime}} = \frac{\alpha(\widetilde{s} - s_c(b_{\ast}))}{2} + \frac{\alpha + 1}{\gamma}.
    \end{equation*}
  \end{enumerate}
\end{lemma}

\begin{proof}
  For convenience, we introduce the following notation:
  \begin{align*}
    & A_1 \coloneq \frac{N-2}{2N} , \quad A_2 \coloneq \frac{s}{N} , \quad A_3 \coloneq \frac{N - 2b_{\ast} + 2\alpha\widetilde{s}}{2N(\alpha + 1)} , \\
    & B_1 \coloneq \frac{1}{2} , \quad B_2 \coloneq \frac{N - 2b_{\ast} + 2\alpha\widetilde{s}+2}{2N(\alpha + 1)} , \\
    & B_3 \coloneq \frac{\widetilde{s}}{N} + \frac{1}{\alpha} , \quad B_4 \coloneq \frac{N-b_{\ast}+\alpha\widetilde{s}}{N(\alpha + 1)}.
  \end{align*} 
  \begin{enumerate}[label=(\roman*)]
    \item We first compare $A_1, A_2, A_3$ with $B_1, B_2$. The inequalities $A_1<B_1$, $A_3<B_2$ are clear, and $A_2<B_1$ follows from $s<\frac{N}{2}$.
    For $A_3 < B_1$, we again use $\widetilde{s}<\frac{N}{2}$ to obtain 
    \begin{align*}
      B_1 - A_3 &= \frac{1}{2} - \frac{N - 2b_{\ast} + 2\alpha\widetilde{s}}{2N(\alpha + 1)} \\
      &= \frac{(N - 2\widetilde{s})\alpha + 2b_{\ast}}{2N(\alpha + 1)} > 0. \\
    \end{align*}
    Next, since $\eta\ll 1$, 
    \begin{align*}
      B_2 - A_1 &= \frac{N - 2b_{\ast} + 2\alpha\widetilde{s}+2}{2N(\alpha + 1)} - \frac{N-2}{2N} \\
      &= \frac{4-2b_{\ast}-(N-2\widetilde{s}-2)\alpha}{2N(\alpha + 1)} \\
      &= \frac{4-2b_{\ast}-(N-2s)\alpha + 2\alpha(1-s+\widetilde{s})}{2N(\alpha + 1)} \\
      &> \frac{\alpha(1-\eta s)}{N(\alpha + 1)} > 0.
    \end{align*}
    Also, since $b_{\ast}<\frac{N}{2}+1-s$ and $\eta$ is sufficiently small, 
    \begin{align*}
      B_2 - A_2 &= \frac{N - 2b_{\ast} + 2\alpha\widetilde{s}+2}{2N(\alpha + 1)} - \frac{s}{N} \\
      &= \frac{N - 2b_{\ast} - 2s + 2 - 2\alpha s\eta}{2N(\alpha + 1)} > 0.
    \end{align*}

    We next consider $B_3$ and $B_4$. If $N\geq 2$, then 
    \begin{equation*}
      B_3 - B_2 = \frac{\alpha(N-2+2b_{\ast}+2\widetilde{s}) + 2N}{2N\alpha(\alpha + 1)} > 0 
    \end{equation*}
    and
    \begin{equation*}
      B_4 - B_2 = \frac{N-2}{2N(\alpha + 1)} \geq 0.
    \end{equation*}
    Thus, a $\rho$ satisfying \refsp{eq:4-1} exists.

    When $N=1$, we have 
    \begin{align*}
      & B_3 - A_2 = \frac{1}{\alpha} - \eta s, \quad B_4 - A_2 = \frac{1-b_{\ast}-s-\eta\alpha s}{\alpha + 1}, \\
      & B_3 - A_3 = \frac{\alpha(1+2b_{\ast}+2(1-\eta)s) + 2}{2\alpha(\alpha + 1)}, \quad B_4 - A_3 = \frac{1}{2(\alpha + 1)}.
    \end{align*}
    All these expressions are positive if $\eta$ is sufficiently small. Hence, a $\rho$ satisfying \refsp{eq:4-1} exists.

    \item Suppose first that $N\geq 3$. Since $A_3<\frac{1}{\rho}<B_2$, we obtain
    \begin{equation*}
      \frac{1}{2} - \frac{b_{\ast}}{N} + \frac{\alpha\widetilde{s}}{N} < \frac{\alpha+1}{\rho} < 1 - \frac{b_{\ast}}{N} + \frac{\alpha\widetilde{s}}{N} - \frac{N-2}{2N}.
    \end{equation*}
    Therefore, $2<r<\frac{2N}{N-2}$. For $N=1,2$, by the inequality $\frac{1}{\rho}<B_4$, we have $2<r<\infty$. Thus $r$ lies in the range of \refsp{\defi:adm}, and by choosing the corresponding exponent $q$, we obtain $(q,r)\in\cA_0$.

    \item By \refsp{eq:4-1}, $r_3 > r_2 > \alpha (>0)$, so $r_2$ and $r_3$ are well-defined. From the definition of $r$ and \refsp{eq:4-2}, 
    \begin{align*}
      \frac{b_{\ast}+s}{N} + \frac{1}{r_{s,b_{\ast}}} = \frac{b_{\ast}}{N} + \frac{1}{r_{0,b_{\ast}}} &= \frac{b_{\ast}}{N} + \alpha\left(\frac{1}{\rho} - \frac{\widetilde{s}}{N} \right) + \frac{1}{\rho} \\
      &= \frac{b_{\ast}}{N} + \frac{\alpha + 1}{\rho} - \frac{\alpha\widetilde{s}}{N} \\
      &= 1 - \frac{1}{r} = \frac{1}{r^{\prime}}.
    \end{align*}
    Since both $(\gamma, \rho)$ and $(q,r)$ are admissible pairs, and
    \begin{equation*}
      s_c(b_\ast) = \frac{N}{2} - \frac{2 - b_{\ast}}{\alpha},
    \end{equation*}
    we also have
    \begin{align*}
      \frac{1}{q^{\prime}} - \frac{\alpha + 1}{\gamma} &= 1 - \frac{N}{4} + \frac{N}{2r} - (\alpha + 1)\left(\frac{N}{4} - \frac{N}{2\rho}\right) \\
      &= 1 - \frac{b_{\ast}}{2} - \frac{\alpha N}{4} + \frac{\alpha \widetilde{s}}{2} \\
      &= \frac{\alpha (\widetilde{s} - s_c(b_{\ast}))}{2}.
    \end{align*}
    The proof is complete. 
  \end{enumerate}
\end{proof}

We now estimate the nonlinear terms using the exponents introduced in \refsp{\lem:4-1}. We first establish the nonlinear estimate for a coefficient satisfying $a\in L^{p_{0,b_{\ast}},\infty}$ and $D^s a\in L^{p_{s,b_{\ast}},\infty}$, where $b_{\ast}\in\{0,b\}$. Let $G(z)\coloneq |z|^{\alpha}z$ for $z\in\C$.

\begin{lemma} \label{lemma:4-2}
  Let $(\gamma,\rho)$ and $(q,r)$ be the admissible pairs given in \refsp{\lem:4-1}. Then, for any finite interval $I\subset\R$ and any $u,v\in L^{\gamma}(I;H_{\rho,2}^{s})$, the following estimates hold:
  \begin{equation} \label{eq:4-3}
    \norm{a|u|^{\alpha}v}_{L^{q^\prime}(I,L^{r^\prime,2})} \lesssim |I|^{\frac{\alpha(\widetilde{s}-s_c(b_{\ast}))}{2}}\norm{u}_{L^{\gamma}(I,H_{\rho,2}^s)}^{\alpha}\norm{v}_{L^{\gamma}(I,L^{\rho,2})},
  \end{equation}
  and
  \begin{equation} \label{eq:4-4}
    \norm{D^s(aG(u))}_{L^{q^\prime}(I,L^{r^\prime,2})} \lesssim |I|^{\frac{\alpha(\widetilde{s}-s_c(b_{\ast}))}{2}}\norm{u}_{L^{\gamma}(I,H_{\rho,2}^{s})}^{\alpha+1}.
  \end{equation}
  Moreover, when $\alpha\geq 1$, 
  \begin{align} 
    & \norm{D^s(a(G(u) - G(v)))}_{L^{q^\prime}(I,L^{r^\prime,2})} \notag \\
    & \hspace{20pt} \lesssim \left[ |I|^{\frac{\widetilde{s}-s_c(b_{\ast})}{2}} \left( \norm{u}_{L^{\gamma}(I,H_{\rho,2}^{s})} + \norm{v}_{L^{\gamma}(I,H_{\rho,2}^{s})} \right) \right]^\alpha \norm{u-v}_{L^{\gamma}(I,H_{\rho,2}^{s})}. \label{eq:4-5}
  \end{align}
  On the other hand, when $0<\alpha<1$, 
  \begin{align}
    &\norm{D^s(a(G(u) - G(v)))}_{L^{q^\prime}(I,L^{r^\prime,2})} \notag\\
    & \hspace{20pt} \lesssim \left[ |I|^{\frac{\widetilde{s}-s_c(b_{\ast})}{2}} \left( \norm{u}_{L^{\gamma}(I,H_{\rho,2}^{s})} + \norm{v}_{L^{\gamma}(I,H_{\rho,2}^{s})} \right) \right]^\alpha \norm{u-v}_{L^{\gamma}(I,H_{\rho,2}^{s})} \notag\\
    & \hspace{40pt} + \left[ |I|^{\frac{\widetilde{s}-s_c(b_{\ast})}{2\eta}} \norm{u-v}_{L^{\gamma}(I,L^{\rho,2})} \right]^{\eta\alpha} \notag\\
    & \hspace{60pt} \times \left( \norm{u}_{L^{\gamma}(I,H_{\rho,2}^{s})} + \norm{v}_{L^{\gamma}(I,H_{\rho,2}^{s})} \right)^{1+(1-\eta)\alpha}. \label{eq:4-6}
  \end{align} 
\end{lemma}

\begin{proof}
  By H\"older's inequality in Lorentz spaces, the Sobolev-Lorentz embedding, and \refsp{\lem:4-1}
  \begin{align*}
    \norm{a|u|^{\alpha}v}_{L^{r^\prime,2}} &\lesssim \norm{a}_{L^{p_{0,b_{\ast}},\infty}} \norm{|u|^{\alpha}v}_{L^{r_{0,b_{\ast}},2}} \\
    &\lesssim \norm{u}_{L^{r_2,\infty}}^{\alpha}\norm{v}_{L^{\rho,2}} \\
    &\lesssim \norm{D^{\widetilde{s}}u}_{L^{\rho,2}}^{\alpha}\norm{v}_{L^{\rho,2}}.
  \end{align*}
  Using the relation between the time exponents in \refsp{\lem:4-1} and H\"older's inequality in time, we obtain
  \begin{align*}
    \norm{a|u|^{\alpha}v}_{L^{q^\prime}(I,L^{r^\prime,2})} &\lesssim |I|^{\frac{\alpha(\widetilde{s} - s_c(b_{\ast}))}{2}} \norm{D^{\widetilde{s}}u}_{L^{\gamma}(I,L^{\rho,2})}^{\alpha}\norm{v}_{L^{\gamma}(I,L^{\rho,2})} \\
    &\lesssim |I|^{\frac{\alpha(\widetilde{s} - s_c(b_{\ast}))}{2}} \norm{u}_{L^{\gamma}(I,H_{\rho,2}^{s})}^{\alpha} \norm{v}_{L^{\gamma}(I,L^{\rho,2})}.
  \end{align*}
  Thus, \refsp{eq:4-3} follows.

  Next, by the fractional Leibniz rule \refsp{\lem:Leibniz}, H\"older's inequality in Lorentz spaces \refsp{\lem:Holder}, and \refsp{eq:4-2}, we have
  \begin{align*}
    &\norm{D^s(aG(u))}_{L^{r^\prime,2}} \\
    &\hspace{30pt} \lesssim \norm{D^s a}_{L^{p_{s,b_{\ast}},\infty}} \norm{G(u)}_{L^{r_{s,b_{\ast}},2}} + \norm{a}_{L^{p_{0,b_{\ast}},\infty}} \norm{D^sG(u)}_{L^{r_{0,b_{\ast}},2}} \\
    &\hspace{30pt} \lesssim \norm{u}_{L^{r_2,\infty}}^{\alpha} \norm{u}_{L^{r_3,2}} + \norm{u}_{L^{r_2,\infty}}^{\alpha} \norm{D^s u}_{L^{\rho,2}} \\
    &\hspace{30pt} \lesssim \norm{D^{\widetilde{s}} u}_{L^{\rho,2}}^{\alpha} \norm{D^s u}_{L^{\rho,2}}.
  \end{align*}
  Here, we used the Sobolev--Lorentz embedding \refsp{\rem:2-8}. For the second term, we applied \refsp{\prop:Fujiwara} with $v=0$ when $0<s<1$. If $s=1$, the classical chain rule and H\"older's inequality yield a similar estimate. Therefore, by H\"older's inequality in time, 
  \begin{align*}
    \norm{D^s(aG(u))}_{L^{q^\prime}(I,L^{r^\prime,2})} &\lesssim |I|^{\frac{\alpha(\widetilde{s} - s_c(b_{\ast}))}{2}} \norm{D^{\widetilde{s}} u}_{L^{\gamma}(I,L^{\rho,2})}^{\alpha} \norm{D^s u}_{L^{\gamma}(I,L^{\rho,2})} \\
    &\lesssim |I|^{\frac{\alpha(\widetilde{s} - s_c(b_{\ast}))}{2}} \norm{u}_{L^{\gamma}(I,H_{\rho,2}^{s})}^{\alpha+1}
  \end{align*}
  which proves \refsp{eq:4-4}.

  We now consider the difference estimate. By the fractional Leibniz rule, together with
  \begin{equation*}
    |G(u) - G(v)| \lesssim (|u|+|v|)^{\alpha}|u-v|
  \end{equation*}
  for the first term and \refsp{\prop:Fujiwara} for the second term, we have
  \begin{align*}
    &\norm{D^s(a(G(u)-G(v)))}_{L^{q^\prime}(I,L^{r^\prime,2})} \\
    & \hspace{50pt} \lesssim \left\|\norm{D^s a}_{L^{p_{s,b_{\ast}},\infty}} \norm{G(u) - G(v)}_{L^{r_{s,b_{\ast}},2}} \right. \\
    & \hspace{100pt} + \left. \norm{a}_{L^{p_{0,b_{\ast}},\infty}} \norm{D^s(G(u) - G(v))}_{L^{r_{0,b_{\ast}},2}}\right\|_{L^{q^\prime}(I)}.
  \end{align*}
  When $\alpha\geq 1$, applying the estimates above and \refsp{\prop:Fujiwara} yields 
  \begin{align*}
    &\norm{D^s(a(G(u)-G(v)))}_{L^{q^\prime}(I,L^{r^\prime,2})} \\
    & \lesssim |I|^{\frac{\alpha(\widetilde{s} - s_c(b_{\ast}))}{2}} \bigg[ \left( \norm{u}_{L^{\gamma}(I,L^{r_2, \infty})} + \norm{v}_{L^{\gamma}(I,L^{r_2, \infty})} \right)^{\alpha} \norm{D^s(u-v)}_{L^{\gamma}(I,L^{\rho, 2})} \\
    & \hspace{20pt} + \left( \norm{D^su}_{L^{\gamma}(I,L^{\rho, 2})} + \norm{D^sv}_{L^{\gamma}(I,L^{\rho, 2})} \right)  \\
    & \hspace{40pt} \times \left( \norm{u}_{L^{\gamma}(I,L^{r_2, \infty})} + \norm{v}_{L^{\gamma}(I,L^{r_2, \infty})} \right)^{\alpha-1} \norm{u-v}_{L^{\gamma}(I,L^{r_2, \infty})} \bigg] \\
    &\lesssim |I|^{\frac{\alpha(\widetilde{s} - s_c(b_{\ast}))}{2}} \bigg[ \left( \norm{D^{\widetilde{s}}u}_{L^{\gamma}(I,L^{\rho, 2})} + \norm{D^{\widetilde{s}}v}_{L^{\gamma}(I,L^{\rho, 2})} \right)^{\alpha} \norm{D^s(u-v)}_{L^{\gamma}(I,L^{\rho, 2})} \\
    & \hspace{20pt} + \left( \norm{D^su}_{L^{\gamma}(I,L^{\rho, 2})} + \norm{D^sv}_{L^{\gamma}(I,L^{\rho, 2})} \right)  \\
    & \hspace{40pt} \times \left( \norm{D^{\widetilde{s}}u}_{L^{\gamma}(I,L^{\rho, 2})} + \norm{D^{\widetilde{s}}v}_{L^{\gamma}(I,L^{\rho, 2})} \right)^{\alpha-1} \norm{D^{\widetilde{s}}(u-v)}_{L^{\gamma}(I,L^{\rho, 2})} \bigg].
  \end{align*}
  Since the interpolation inequality $\norm{D^{\widetilde{s}}f} \lesssim \norm{D^sf}^{1-\eta}\norm{f}^{\eta}$ holds, \refsp{eq:4-5} follows. Suppose next that $0<\alpha <1$ and $s<1$. Taking $\theta = \alpha$ in \refsp{\prop:Fujiwara} and using $\widetilde{s} = (1-\eta)s$, we get
  \begin{align*}
    & \norm{D^s(a(G(u) - G(v)))}_{L^{q^\prime}(I,L^{r^\prime,2})} \\
    & \hspace{20pt} \lesssim |I|^{\frac{\alpha(\widetilde{s} - s_c(b_{\ast}))}{2}} \left( \norm{D^{\widetilde{s}}u}_{L^{\gamma}(I,L^{\rho, 2})} + \norm{D^{\widetilde{s}}v}_{L^{\gamma}(I,L^{\rho, 2})} \right)^{\alpha} \\
    & \hspace{200pt} \times \norm{D^s(u-v)}_{L^{\gamma}(I,L^{\rho, 2})} \\
    & \hspace{40pt} + |I|^{\frac{\alpha(\widetilde{s} - s_c(b_{\ast}))}{2}} \left( \norm{D^su}_{L^{\gamma}(I,L^{\rho, 2})} + \norm{D^sv}_{L^{\gamma}(I,L^{\rho, 2})} \right) \\
    & \hspace{200pt} \times  \norm{D^{\widetilde{s}}(u-v)}_{L^{\gamma}(I,L^{\rho, 2})}^{\alpha} \\
    & \hspace{20pt} \lesssim |I|^{\frac{\alpha(\widetilde{s} - s_c(b_{\ast}))}{2}} \left( \norm{u}_{L^{\gamma}(I,H_{\rho, 2}^{s})} + \norm{v}_{L^{\gamma}(I,H_{\rho, 2}^{s})} \right)^{\alpha} \norm{D^s(u-v)}_{L^{\gamma}(I,L^{\rho, 2})} \\
    & \hspace{40pt} + |I|^{\frac{\alpha(\widetilde{s} - s_c(b_{\ast}))}{2}} \left( \norm{u}_{L^{\gamma}(I,H_{\rho, 2}^{s})} + \norm{v}_{L^{\gamma}(I,H_{\rho, 2}^{s})} \right) \\
    & \hspace{60pt} \times \norm{D^{s}(u-v)}_{L^{\gamma}(I,L^{\rho, 2})}^{(1-\eta)\alpha} \norm{u-v}_{L^{\gamma}(I,L^{\rho, 2})}^{\eta\alpha} \\
    & \hspace{20pt} \lesssim |I|^{\frac{\alpha(\widetilde{s} - s_c(b_{\ast}))}{2}} \left( \norm{u}_{L^{\gamma}(I,H_{\rho, 2}^{s})} + \norm{v}_{L^{\gamma}(I,H_{\rho, 2}^{s})} \right)^{\alpha} \norm{u-v}_{L^{\gamma}(I,H_{\rho, 2}^{s})} \\
    & \hspace{40pt} + |I|^{\frac{\alpha(\widetilde{s} - s_c(b_{\ast}))}{2}} \left( \norm{u}_{L^{\gamma}(I,H_{\rho, 2}^{s})} + \norm{v}_{L^{\gamma}(I,H_{\rho, 2}^{s})} \right) \\
    & \hspace{60pt} \times \norm{D^{s}(u-v)}_{L^{\gamma}(I,L^{\rho, 2})}^{(1-\eta)\alpha} \norm{u-v}_{L^{\gamma}(I,L^{\rho, 2})}^{\eta\alpha} \\
    & \hspace{20pt} \lesssim |I|^{\frac{\alpha(\widetilde{s} - s_c(b_{\ast}))}{2}} \left( \norm{u}_{L^{\gamma}(I,H_{\rho, 2}^{s})} + \norm{v}_{L^{\gamma}(I,H_{\rho, 2}^{s})} \right)^{\alpha} \norm{u-v}_{L^{\gamma}(I,H_{\rho, 2}^{s})} \\
    & \hspace{40pt} + |I|^{\frac{\alpha(\widetilde{s} - s_c(b_{\ast}))}{2}} \left( \norm{u}_{L^{\gamma}(I,H_{\rho, 2}^{s})} + \norm{v}_{L^{\gamma}(I,H_{\rho, 2}^{s})} \right)^{1+(1-\eta)\alpha}  \\
    & \hspace{200pt} \times \norm{u-v}_{L^{\gamma}(I,L^{\rho, 2})}^{\eta\alpha}.
  \end{align*} 
  This is exactly \refsp{eq:4-6}. If $s=1$, the same estimate follows from the classical chain rule. This completes the proof of \refsp{\lem:4-2}. 
\end{proof}

\begin{proof}[Proof of \refsp{\thm:1-1}]
  We first prove the existence of solutions by using the mapping $\Phi$ defined in \refsp{eq:1-5}. Following \refsp{\rem:2-15}, we estimate the bounded and singular terms separately. For each $b_{\ast}\in\{0,b\}$, let $(\gamma_{b_{\ast}},\rho_{b_{\ast}})$ and $(q_{b_{\ast}},r_{b_{\ast}})$ be the admissible pairs given by \refsp{\lem:4-1}. We also set
  \begin{equation*}
    \beta_{b_{\ast}} \coloneq \frac{\alpha(\widetilde{s}-s_c(b_{\ast}))}{2}, \quad \beta \coloneq \frac{\alpha(\widetilde{s}-s_c)}{2}.
  \end{equation*}
  Since $s_c(b_{\ast})\leq s_c$, we have $T^{\beta_{b_{\ast}}}\leq T^{\beta}$ for $0<T\leq 1$.

  Define $X(M,T)$ as the set of all functions $u$ satisfying
  \begin{equation*}
    u \in L^{\infty}(0,T;H^{s}) \cap \bigcap_{b_{\ast}\in\{0,b\}} L^{\gamma_{b_{\ast}}}(0,T;H_{\rho_{b_{\ast}},2}^{s})
  \end{equation*}
  together with
  \begin{equation*}
    \norm{u}_{L_T^\infty H^s} + \sum_{b_{\ast}\in\{0,b\}} \norm{u}_{L_T^{\gamma_{b_{\ast}}} H_{\rho_{b_{\ast}},2}^{s}} \leq M
  \end{equation*}
  and, for $u,v\in X(M,T)$,
  \begin{equation*}
    d(u,v) \coloneq \norm{u-v}_{L_T^\infty L^2} + \sum_{b_{\ast}\in\{0,b\}} \norm{u-v}_{L_T^{\gamma_{b_{\ast}}} L^{\rho_{b_{\ast}},2}}.
  \end{equation*}
  { By the reflexivity of the Sobolev--Lorentz spaces and the weak-* compactness of bounded sets in $L^\infty(0,T;H^s)$, we can show that $(X(M,T),d)$ is a complete metric space by the same argument as in \cite[Proof of Theorem 1.2]{A-T}.}
  
  In the following argument, we use the decomposition in \refsp{\rem:2-15} and estimate the bounded and singular terms separately. The estimates corresponding to $b_{\ast} = 0$ and  $b_{\ast} = b$ contain the time fafctors $T^{\beta_0}$ and $T^{\beta_b}$, respectively. Since both are bounded by $T^{\beta}$ for $0 < T \leq 1$, the same contraction argument applies to the general coefficient.
  
  Assume that $\hsnorm{u_0}\leq R$. By the Strichartz estimates and \refsp{\lem:4-2},
  \begin{align*}
    &\norm{\Phi(u)}_{L_T^\infty H^s} + \norm{\Phi(u)}_{L_T^\gamma H_{\rho,2}^s} \\
    &\hspace{30pt} \lesssim \hsnorm{u_0} + \norm{aG(u)}_{L_T^{q^\prime}L^{r^\prime,2}} + \norm{D^s(aG(u))}_{L_T^{q^\prime}L^{r^\prime,2}} \\
    &\hspace{30pt} \lesssim R + T^{\beta} \norm{u}_{L_T^\gamma H_{\rho,2}^s}^{\alpha+1}.
  \end{align*}
  Hence, if $u\in X(M,T)$,
  \begin{equation*}
    \norm{\Phi(u)}_{L_T^\infty H^s} + \norm{\Phi(u)}_{L_T^\gamma H_{\rho,2}^s} \leq CR + CT^{\beta} M^{\alpha+1},
  \end{equation*}
  where $C>0$ is a constant independent of $u_0$, $R$, $M$ and $T$. For the contraction estimate, we use the inequality
  \begin{equation*}
    |G(u)-G(v)|\lesssim (|u|^\alpha+|v|^\alpha)|u-v|.
  \end{equation*}
  By the Strichartz estimates, the above inequality, and (4.3) , we obtain, for $u, v \in X(M,T)$, 
  \begin{align*}
    d(\Phi(u),\Phi(v)) &= \norm{\Phi(u) - \Phi(v)}_{L_T^\infty L^2} + \norm{\Phi(u) - \Phi(v)}_{L_T^\gamma L^{\rho,2}} \\
    &\lesssim \norm{a(G(u) - G(v))}_{L_T^{q^\prime}L^{r^\prime,2}} \\
    &\lesssim T^{\beta} \left( \norm{u}_{L_T^\gamma H_{\rho,2}^s}^{\alpha} + \norm{v}_{L_T^\gamma H_{\rho,2}^s}^{\alpha} \right) \norm{u-v}_{L_T^\gamma L^{\rho,2}} \\
    &\lesssim 2 T^{\beta} M^{\alpha} d(u,v).
  \end{align*}
  Now take $M= 2CR$ and choose $T=T(R)>0$ sufficiently small so that $CT^{\beta} M^{\alpha}\leq \frac{1}{4}$. Then $\Phi$ maps $X(M,T)$ into itself and is a contraction on $X(M,T)$. Therefore, the Banach fixed-point theorem gives a unique solution $u\in X(M,T)$ to \refsp{eq:1-1}. The time continuity of the Duhamel term also yields $u\in C([0,T];H^s)$. By repeating the same argument with a shifted initial time, this local solution can be extended up to the maximal existence time $T_{\max}(u_0)>0$. 

  \begin{enumerate}[label=(\roman*)]
    \item Fix any $T\in[0,T_{\max}(u_0))$, and let 
    \begin{equation*}
      u,v\in C([0,T];H^{s}(\R^N)) \cap \bigcap_{b_{\ast}\in\cB} L^{\gamma_{b_{\ast}}}(0,T;H_{\rho_{b_{\ast}},2}^{s})
    \end{equation*}
    be two solutions of \refsp{eq:1-1} with the same initial data $u_0$. As before, we sum the estimates over all coefficient components and write only the exponents corresponding to a single component in order to simplify the notation. We divide $[0,T]$ into finitely many intervals $I_j \coloneq [t_{j-1}, t_j]$, $0 = t_0 < t_1 < \cdots < t_m = T$ such that $t_j - t_{j-1} < h$ for every $j$. For $j=1$, since $u(t_0) = u(0) = v(0) = v(t_0)$, the same nonlinear estimate as above implies 
    \begin{align*}
      &\norm{u - v}_{L^{\infty}(I_1;L^{2})} + \norm{u - v}_{L^{\gamma}(I_1;L^{\rho,2})} \\
      & \hspace{30pt} \lesssim |I_1|^{\beta}\left( \norm{u}_{L^{\gamma}(I_1;H_{\rho,2}^{s})}^{\alpha} + \norm{v}_{L^{\gamma}(I_1;H_{\rho,2}^{s})}^{\alpha} \right) \norm{u - v}_{L^{\gamma}(I_1;L^{\rho,2})} \\
      & \hspace{30pt} \lesssim h^{\beta} \left( \norm{u}_{L^{\gamma}(I_1;H_{\rho,2}^{s})}^{\alpha} + \norm{v}_{L^{\gamma}(I_1;H_{\rho,2}^{s})}^{\alpha} \right) \norm{u - v}_{L^{\gamma}(I_1;L^{\rho,2})}.
    \end{align*}
    Since $\beta>0$, we may choose $h>0$ sufficiently small so that 
    \begin{equation*}
      C h^\beta\left( \norm{u}_{L_T^\gamma H_{\rho,2}^s}^{\alpha} + \norm{v}_{L_T^\gamma H_{\rho,2}^s}^{\alpha} \right)<\frac{1}{2}.
    \end{equation*}
    It then follows that $u \equiv v\ \aeae$ on $I_1$. In particular, $u(t_1) = v(t_1)$, so the same argument can be applied on $I_2$. Repeating this procedure over all the intervals $I_j$, we obtain
    \begin{equation*}
      u \equiv v \hspace{10pt} \aeae \hspace{5pt} \mathrm{on} \hspace{5pt} [0,T].
    \end{equation*}
    This proves uniqueness.
    
    \item Let $(\hat{q},\hat{r})$ be any admissible pair. By the Strichartz estimates and the nonlinear estimate above,
    \begin{equation*}
      \norm{u}_{L_T^{\hat{q}}H_{\hat{r},2}^{s}} \lesssim \hsnorm{u_0} + T^{\beta} \norm{u}_{L_T^{\gamma}H_{\rho,2}^s}^{\alpha + 1} < \infty.
    \end{equation*}
    Therefore, $u\in L^{\hat{q}}(0,T;H_{\hat{r},2}^{s})$ for every admissible pair $(\hat{q},\hat{r})$. 
    
    \item Let $R>0$ and suppose that $\hsnorm{u_0}, \hsnorm{v_0}\leq R$. By the Duhamel formula \refsp{eq:1-5},
    \begin{equation*}
      u(t) - v(t) = e^{it\Delta}(u_0 - v_0) - i\int_0^t e^{i(t-\tau)\Delta}a(x)(G(u) - G(v))(\tau) d\tau.
    \end{equation*}
    As in the contraction mapping argument, we obtain
    \begin{align*}
      &\norm{u-v}_{L_T^{\infty}L^2} + \norm{u-v}_{L_T^{\gamma}L^{\rho,2}} \\
      & \hspace{30pt} \lesssim \lnorm{u_0 - v_0} + T^{\beta} \left( \norm{u}_{L_T^{\gamma}H_{\rho,2}^{s}} + \norm{v}_{L_T^{\gamma}H_{\rho,2}^{s}} \right)^{\alpha} \norm{u-v}_{L_T^{\gamma}L^{\rho,2}}. 
    \end{align*}
    As in the local existence argument, we choose $T>0$ sufficiently small so that the last term on the right-hand side can be absorbed into the left-hand side. Hence, 
    \begin{equation*}
      \norm{u-v}_{L_T^{\infty}L^2} + \norm{u-v}_{L_T^{\gamma}L^{\rho,2}} \lesssim \lnorm{u_0 - v_0}.
    \end{equation*}

    When $\alpha\geq 1$, the Duhamel formula, Strichartz estimates, and \refsp{\lem:4-2} yield
    \begin{align*}
      &\sup_{0\leq t\leq T}\norm{u(t)-v(t)}_{H^{s}} + \norm{u-v}_{L_T^{\gamma}H_{\rho,2}^s} \\
      & \hspace{30pt} \lesssim \hsnorm{u_0 - v_0} + T^{\beta} \left( \norm{u}_{L_T^{\gamma}H_{\rho,2}^{s}} + \norm{v}_{L_T^{\gamma}H_{\rho,2}^{s}} \right)^\alpha \norm{u-v}_{L_T^{\gamma}H_{\rho,2}^{s}}.
    \end{align*}
    Taking $T>0$ smaller if necessary, we can again absorb the last term into the left-hand side, and thus 
    \begin{equation*}
      \sup_{0\leq t\leq T}\norm{u(t)-v(t)}_{H^{s}} + \norm{u-v}_{L_T^{\gamma}H_{\rho,2}^s} \lesssim \hsnorm{u_0 - v_0}.
    \end{equation*}

    Finally, suppose that $0<\alpha<1$. Applying \refsp{\lem:4-2} in a similar way, we obtain
    \begin{align*}
      & \sup_{0\leq t\leq T}\norm{D^s(u(t)-v(t))}_{L^{2}} + \norm{D^s(u-v)}_{L_T^{\gamma}L^{\rho,2}} \\
      & \hspace{20pt} \lesssim \lnorm{D^s(u_0 - v_0)} + T^{\beta} \left( \norm{u}_{L_T^{\gamma}H_{\rho,2}^{s}} + \norm{v}_{L_T^{\gamma}H_{\rho,2}^{s}} \right)^\alpha \norm{u-v}_{L_T^{\gamma}H_{\rho,2}^{s}} \\
      & \hspace{40pt} + T^{\beta} \norm{u-v}_{L_T^{\gamma}L^{\rho,2}}^{\eta\alpha} \left( \norm{u}_{L_T^{\gamma}H_{\rho,2}^{s}} + \norm{v}_{L_T^{\gamma}H_{\rho,2}^{s}} \right)^{1+(1-\eta)\alpha}
    \end{align*}
    The term containing $D^s(u-v)$ is included in the $H^s$-norm on the right-hand side. Hence, by taking $T>0$ sufficiently small, we get
    \begin{align*}
      &\sup_{0\leq t\leq T}\norm{u(t)-v(t)}_{H^{s}} + \norm{u-v}_{L_T^{\gamma}H_{\rho,2}^s} \lesssim \hsnorm{u_0 - v_0} + C_R \hsnorm{u_0 - v_0}^{\eta\alpha}.
    \end{align*}
    
    In both cases for $\alpha$, the same argument as in (ii) allows us to obtain the corresponding estimate for any admissible pair. Therefore, \refsp{eq:6} follows.
  \end{enumerate}
\end{proof}

\section{Local well-posedness in the critical case} \label{sec:5}

In this section, we prove \refsp{\thm:1-2} in the critical case
\begin{equation*}
  \alpha = \alpha_c \coloneq \frac{4-2b}{N-2s}.
\end{equation*}
Throughout this section, we assume the hypotheses of \refsp{\thm:1-2} and consider the case $0<s\leq1$. The case $s=0$ can be treated in the same way by omitting the estimates involving derivatives. When $b > 0$, we have
\begin{equation*}
  \alpha_c < \frac{4}{N-2s}.
\end{equation*}
Thus, the terms corresponding to $b_{\ast} = 0$ can be treated by the subcritical estimates in \refsp{sec:4}. Therefore, we focus on the term corresponding to $b_{\ast} = b$ below. When $b=0$, no separate subcritical terms occur.

We shall choose the auxiliary parameter $0<s_0<\min\{\alpha s,s\}$ later. 

We consider the following condition on:
\begin{align} \label{eq:5-1}
  & \max\left\{ \frac{N-2}{2N}, \frac{s}{N}, \frac{N - 2b + 2(\alpha+1)s - 2s_0}{2N(\alpha + 1)} \right\} < \frac{1}{\rho} \notag\\
  & \hspace{120pt}  < \min\Bigg\{ \frac{1}{2} - \frac{1}{N(\alpha + 1)}, \frac{N - b + \alpha s}{N(\alpha + 1)} \Bigg\}.
\end{align}

\begin{lemma} \label{lemma:5-1}
  The following statements hold.
  \begin{enumerate}[label=(\roman*)]
    \item There exists an admissible pair $(\gamma, \rho)$ satisfying \refsp{eq:5-1}. 
    \item For any such pair $(\gamma, \rho)$, define
    \begin{equation*}
      \frac{1}{r} = 1 - \frac{b}{N} - \frac{\alpha + 1}{\rho} + \frac{\alpha s}{N}.
    \end{equation*} 
    Then, there exists $q$ such that $(q,r)\in\cA_0$. 
    \item For the pairs $(\gamma, \rho)$ and $(q,r)$ above, we have
    \begin{equation*}
      \frac{1}{q} = 1 - \frac{\alpha+1}{\gamma}.
    \end{equation*}
  \end{enumerate}
\end{lemma}

\begin{proof}
  For convenience, we use the notation
  \begin{align*}
    & A_4 \coloneq \frac{N-2}{2N} , \quad A_5 \coloneq \frac{s}{N} , \quad A_6 \coloneq \frac{N - 2b + 2(\alpha + 1)s - 2s_0}{2N(\alpha + 1)} , \\
    & B_5 \coloneq \frac{1}{2} - \frac{1}{N(\alpha + 1)} , \quad B_6 \coloneq \frac{N - b + \alpha s}{N(\alpha + 1)}
  \end{align*}
  when necessary. 
  \begin{enumerate}[label=(\roman*)]
    \item It is enough to compare $A_4, A_5, A_6$ with $B_5, B_6$. 
    \begin{enumerate}[label=(\alph*)]
      \item For $A_4 < B_5$, we have
      \begin{align*}
        B_5 - A_4 &= \frac{N(\alpha + 1) - 2 - (\alpha+1)(N-2)}{2N(\alpha + 1)} \\
        &= \frac{\alpha}{N(\alpha + 1)} > 0.
      \end{align*}
      \item For $A_5 < B_5$, since $b<\frac{N}{2}+1-s$, we have
      \begin{align*}
        B_5 - A_5 &=  \frac{(N - 2s)(\alpha + 1) - 2}{2N(\alpha + 1)} \\
        &= \frac{2-2b+N-2s}{2N(\alpha + 1)} > 0.
      \end{align*}
      \item For $A_6 < B_5$, since $0<s\leq 1$ and $0 < s_0$, we have
      \begin{align*}
        B_5 - A_6 &= \frac{(N - 2s)\alpha - 2 + 2b - 2s + 2s_0}{2N(\alpha + 1)} \\
        &= \frac{1-s+s_0}{N(\alpha + 1)} > 0.
      \end{align*}
      \item For $A_4 < B_6$, we note that $\alpha>1$ when $N=1$. Then
      \begin{align*}
        B_6 - A_4 &= \frac{N-2b+2-(N-2s-2)\alpha}{2N(\alpha + 1)} \\
        &= \frac{N-2+2\alpha}{2N(\alpha + 1)} > 0.
      \end{align*}
      \item For $A_5 < B_6$, since $b<N-s$, we have
      \begin{align*}
        B_6 - A_5 &= \frac{N-b+\alpha s-s(\alpha+1)}{N(\alpha + 1)} \\
        &= \frac{N-b-s}{N(\alpha + 1)} > 0.
      \end{align*}
      \item For $A_6 < B_6$, since $s<\frac{N}{2}$ and $s_0>0$, we have
      \begin{align*}
        B_6 - A_6&= \frac{N-b+\alpha s}{N(\alpha + 1)} - \frac{N - 2b + 2(\alpha+1)s - 2s_0}{2N(\alpha + 1)} \\
        &= \frac{N-2s+2s_0}{2N(\alpha + 1)} > 0.
      \end{align*}
    \end{enumerate}
    Therefore, we can choose $\gamma$ such that $(\gamma, \rho)$ is an admissible pair.

    \item By the condition on $A_6$, we have 
    \begin{align*}
      \frac{1}{r} &< 1 - \frac{b}{N} - \frac{N - 2b + 2(\alpha + 1)s - 2s_0}{2N} + \frac{\alpha s}{N} \\
      &= \frac{1}{2} - \frac{s - s_0}{N} < \frac{1}{2}.
    \end{align*}
    On the other hand, by the condition on $B_5$, 
    \begin{align*}
      \frac{1}{r} &> 1 - \frac{b}{N} - (\alpha + 1)\left( \frac{1}{2} - \frac{1}{N(\alpha + 1)} \right) + \frac{\alpha s}{N} \\
      & = \frac{N-2}{2N}.
    \end{align*}
    Moreover, the condition on $B_6$ gives $r>0$. Therefore, 
    \begin{equation*}
      2 < r < \begin{dcases}
        \frac{2N}{N-2} & N\geq 3, \\
        \infty & N=1,2.
      \end{dcases}
    \end{equation*}
    Hence, by \refsp{\defi:adm}, we can choose $q$ such that $(q,r)\in\cA_0$. 

    \item Since $(q,r)$ and $(\gamma,\rho)$ are admissible pairs, we have
    \begin{align*}
      \frac{1}{q} &= \frac{N}{2} \left( \frac{1}{2} - \frac{1}{r} \right) \\
      &= \frac{N}{2} \left( \frac{b}{N} + \frac{\alpha + 1}{\rho} - \frac{1}{2} - \frac{\alpha s}{N} \right) \\
      &= \frac{2b+\alpha(N-2s)}{4} - \frac{\alpha+1}{\gamma} \\
      &= 1 - \frac{\alpha+1}{\gamma} .
    \end{align*}
  \end{enumerate}
\end{proof}

For the nonlinear estimates in the critical case, we use a $\dot{H}^{s-s_0}$-admissible pair $(\gamma_0,\rho_0)$. We next show that such exponents can indeed be chosen. 

\begin{lemma} \label{lemma:5-2}
  Let $0<\alpha < 1$. Define $(\gamma_2,\rho_2)$ by
  \begin{equation*}
    \gamma_2 \coloneq \gamma, \quad \frac{1}{\rho_2} \coloneq \frac{1}{\rho}-\frac{s-s_0}{N}.
  \end{equation*}
  Then $(\gamma_2,\rho_2)$ is a $\dot{H}^{s-s_0}$-admissible pair. Moreover, for a suitable $\theta\in(0,1)$, the pairs $(\gamma_1,\rho_1)$, $(\gamma_1,\rho_1)$, and $(\widetilde{\gamma}_0, \widetilde{\rho}_0)$ defined in the proof below satisfy
  \begin{equation*}
    (\gamma_1,\rho_1) \in\cA_0, \quad (\gamma_0,\rho_0) \in\cA_{s-s_0}, 
  \end{equation*}
  and $(\widetilde{\gamma}_0, \widetilde{\rho}_0)$ is an acceptable pair. Furthermore, $(\gamma_0, \rho_0)$ and $(\widetilde{\gamma}_0, \widetilde{\rho}_0)$ satisfy the assumptions of \refsp{\prop:2-12}. 
\end{lemma}

\begin{proof}
  When $N=1$, we have $\alpha>1$. Thus, we only consider $N\geq 2$ below.

  Choose $\theta\in(0,1)$ such that
  \begin{equation*}
     \frac{1-(s-s_0)}{N(\alpha+1)\left(\frac{1}{2}-\frac{1}{\rho}\right)} < \theta < \min\left\{1, \frac{1-(s-s_0)}{N\left(\frac{1}{2}-\frac{1}{\rho}\right)}
  \right\}.
  \end{equation*}
  By the condition \refsp{eq:5-1} and $\alpha + 1>1$, the interval above is non-empty, so such a $\theta$ can be chosen. Define $(\gamma_1,\rho_1)$ by
  \begin{equation*}
    \frac{1}{\gamma_1} \coloneq \frac{\theta}{\gamma}, \quad \frac{1}{\rho_1} \coloneq \frac{1 - \theta}{2} + \frac{\theta}{\rho}.
  \end{equation*}
  We further define $(\gamma_0,\rho_0)$ and $(\widetilde{\gamma}_0, \widetilde{\rho}_0)$ by 
  \begin{gather*}
    \gamma_0:=\gamma_1, \quad \frac{1}{\rho_0} := \frac{1}{\rho_1} -\frac{s-s_0}{N} \\
    \frac{1}{\widetilde{\gamma}_0} := 1 - \frac{\alpha}{\gamma_1} - \frac{1}{\gamma_0}, \quad 1-\frac{1}{\widetilde{\rho}_0} := \frac{b}{N} + \alpha\left( \frac{1}{\rho_1} - \frac{s}{N} \right) + \frac{1}{\rho_0}.
  \end{gather*} 

  We first verify that $(\gamma_1,\rho_1)\in\cA_0$. Since $(\gamma,\rho)$ is an admissible pair,
  \begin{align*}
    \frac{2}{\gamma_1}+\frac{N}{\rho_1} &= \frac{2\theta}{\gamma} + \frac{N(1-\theta)}{2} + \frac{N\theta}{\rho} \\
    &= \theta\left( \frac{2}{\gamma}+\frac{N}{\rho} \right) + (1-\theta)\frac{N}{2} = \frac{N}{2}.
  \end{align*}
  Also, since $(\gamma,\rho)\in\cA_0$,
  \begin{equation*}
    0<\frac{1}{2}-\frac{1}{\rho}<\frac{1}{N}.
  \end{equation*}
  Therefore,
  \begin{equation*}
    \frac{N-2}{2N}<\frac{1}{\rho_1}=\frac{1}{2}-\theta\left( \frac{1}{2}-\frac{1}{\rho} \right) < \frac{1}{2}.
  \end{equation*}
  Hence, $(\gamma_1,\rho_1)\in\cA_0$. 

  Next, we verify that $(\gamma_0,\rho_0)\in\cA_{s-s_0}$. By definition,
  \begin{equation*}
    \frac{1}{\rho_0} = \frac{1}{2} - \theta\left( \frac{1}{2}-\frac{1}{\rho} \right) - \frac{s-s_0}{N}. 
  \end{equation*}
  Thus, by the choice of $\theta$,
  \begin{equation*}
    \frac{1}{2}-\frac{1}{N} < \frac{1}{\rho_0} < \frac{1}{2}-\frac{s-s_0}{N}.
  \end{equation*}
  Moreover, since $\gamma_0=\gamma_1$ and $(\gamma_1,\rho_1)\in\cA_0$,
  \begin{align*}
    \frac{2}{\gamma_0} +\frac{N}{\rho_0} &= \frac{2}{\gamma_1} + \frac{N}{\rho_1} - (s-s_0) \\
    &= \frac{N}{2} - (s-s_0).
  \end{align*}
  Therefore, $(\gamma_0,\rho_0)\in\cA_{s-s_0}$.

  On the other hand, for $(\widetilde{\gamma}_0,\widetilde{\rho}_0)$, using the critical relation $\alpha(N-2s)=4-2b$, we obtain
  \begin{align*}
    \frac{1}{\widetilde{\gamma}_0} &= 1- \frac{N(\alpha+1)\theta}{2} \left( \frac{1}{2} - \frac{1}{\rho} \right) \\
    \frac{1}{\widetilde{\rho}_0} &= \frac{1}{2} - \frac{2-(s-s_0)}{N} + (\alpha + 1) \theta \left( \frac{1}{2} - \frac{1}{\rho} \right).
  \end{align*}
  By the choice of $\theta$ and $0<\alpha<1$, we have
  \begin{equation*}
    0 < \frac{1}{\widetilde{\rho}_0} < \frac{1}{2}, \quad 0 < \frac{1}{\widetilde{\gamma}_0} < 1.
  \end{equation*}
  Moreover, 
  \begin{equation*}
    \frac{N}{2} - \frac{N}{\widetilde{\rho}_0} - \frac{1}{\widetilde{\gamma}_0} = 1-(s-s_0) - \frac{N(\alpha+1)\theta}{2} \left( \frac{1}{2} - \frac{1}{\rho} \right) >0.
  \end{equation*}
  Hence, $(\widetilde{\gamma}_0,\widetilde{\rho}_0)$ is an acceptable pair.

  We next check that $(\gamma_0,\rho_0)$, $(\widetilde{\gamma}_0,\widetilde{\rho}_0)$ satisfy the assumptions of \refsp{\prop:2-12}. By the definitions of the indices,
  \begin{align*}
    \frac{1}{\gamma_0} +\frac{1}{\widetilde{\gamma}_0} = 1-\frac{\alpha}{\gamma_1} < 1 \\
    \frac{1}{\gamma_0} + \frac{1}{\widetilde{\gamma}_0} = \frac{N}{2} \left( 1 - \frac{1}{\rho_0} - \frac{1}{\widetilde{\rho}_0} \right).
  \end{align*}
  Also,
  \begin{align*}
    \frac{N}{\widetilde{\rho}_0} - \frac{N-2}{\rho_0} &= \left( N(\alpha + 2) - 2 \right) \theta \left( \frac{1}{2} - \frac{1}{\rho} \right) + \frac{2(N-1)(s-s_0)}{N} - 1 \\
    &= \left( N(\alpha + 2) - 2 \right) \left\{ \theta \left( \frac{1}{2} - \frac{1}{\rho} \right) - \frac{1 - (s - s_0)}{N(\alpha + 1)} \right\} \\
    & \hspace{100pt} + \frac{(N - 2)\{ 1+\alpha(s - s_0) \}}{N(\alpha + 1)} > 0 \\
    \frac{N}{\rho_0} - \frac{N-2}{\widetilde{\rho}_0} &= 1 + \frac{2(N-2)}{N} - \frac{\{N + (N-2)\}(s - s_0)}{N} \\
    & \hspace{100pt} - \{N + (\alpha + 1)(N - 2)\}\theta\left(\frac{1}{2} - \frac{1}{\rho}\right) \\
    & = \{ (N-2)\alpha + 2(N-1) \} \left\{ \frac{1-(s-s_0)}{N} - \theta\left(\frac{1}{2} - \frac{1}{\rho}\right) \right\} \\
    & \hspace{100pt} + \frac{(N-2)\{1 - \alpha\{1 - (s - s_0)\}\}}{N} > 0.
  \end{align*}
  Therefore, $(\gamma_0,\rho_0)$ and $(\widetilde{\gamma}_0,\widetilde{\rho}_0)$ satisfy the assumptions of \refsp{\prop:2-12}.

  { It remains to verify that $(\gamma_2,\rho_2)$ is a $\dot H^{s-s_0}$-admissible pair.} By \refsp{eq:5-1}, 
  \begin{equation*}
    0 < \frac{s_0}{N} < \frac{1}{\rho_2} < \frac{1}{2}.
  \end{equation*}
  Hence, $\rho_2$ is well-defined. Also,
  \begin{equation*}
    \frac{2}{\gamma_2} + \frac{N}{\rho_2} = \frac{2}{\gamma}+\frac{N}{\rho}-(s-s_0) = \frac{N}{2}-(s-s_0).
  \end{equation*}

  Furthermore, by \refsp{\lem:5-1}, 
  \begin{align*}
    \frac{1}{q'} &= \frac{\alpha+1}{\gamma} = \frac{\alpha}{\gamma_2} + \frac{1}{\gamma} \\
    \frac{1}{r'} &= \frac{b}{N} + \frac{\alpha+1}{\rho} - \frac{\alpha s}{N} \\
    &= \frac{b}{N} + \alpha\left( \frac{1}{\rho_2} - \frac{s_0}{N} \right) + \frac{1}{\rho}. 
  \end{align*}
  Thus, we obtain \refsp{\lem:5-2}.
\end{proof}

\begin{lemma} \label{lemma:5-3}
  Let $0<\alpha<1$. Suppose that $a\in L^{p_{0,b},\infty}$ and $D^sa\in L^{p_{s,b},\infty}$. Let $(\gamma,\rho)$ and $(q,r)$ be the admissible pairs given in \refsp{\lem:5-1}, and let $(\gamma_0,\rho_0)$, $(\gamma_1,\rho_1)$, $(\gamma_2,\rho_2)$, and $(\widetilde{\gamma}_0, \widetilde{\rho}_0)$ be the exponents given in \refsp{\lem:5-2}. Then, for any finite interval $I\subset\R$, 
  \begin{align}
    &\norm{ D^{s_0} \left( a(x)|u|^{\alpha}w \right) }_{L^{\widetilde{\gamma}_0^{\prime}}(I,L^{\widetilde{\rho}_0^{\prime},2})} \lesssim \norm{u}_{L^{\gamma_1}(I,H_{\rho_1,2}^{s})}^{\alpha} \norm{D^{s_0}w}_{L^{\gamma_0}(I,L^{\rho_0,2})} \label{eq:5-2} \\
    &\norm{D^s\left[ a(x)\left( |u|^{\alpha}u - |v|^{\alpha}v \right) \right]}_{L^{q^\prime}(I,L^{r^\prime,2})} \notag \\
    & \hspace{30pt} \lesssim \left( \norm{u}_{L^{\gamma}(I,H_{\rho,2}^{s})} + \norm{v}_{L^{\gamma}(I,H_{\rho,2}^{s})} \right)^{\alpha} \norm{u - v}_{L^{\gamma}(I,H_{\rho,2}^{s})} \notag \\
    & \hspace{60pt} + \norm{D^{s_0}(u - v)}_{L^{\gamma_2}(I,L^{\rho_2,2})}^{\alpha} \left( \norm{u}_{L^{\gamma}(I,H_{\rho,2}^{s})} + \norm{v}_{L^{\gamma}(I,H_{\rho,2}^{s})} \right). \label{eq:5-3}
  \end{align}
\end{lemma}

\begin{proof}
  We first prove \refsp{eq:5-2}. Fix $\sigma$ such that $\frac{s_0}{\alpha} < \sigma < s$.  Also, set
  \begin{equation*}
    \frac{1}{r_4} = \frac{1}{\rho_1} - \frac{s}{N} = \frac{1}{\rho_0} - \frac{s_0}{N}, \quad \frac{1}{r_5} = \frac{\alpha}{r_4} + \frac{s_0}{N}, \quad \frac{1}{r_6} = \frac{1}{\rho_1} - \frac{s - \sigma}{N}.
  \end{equation*}
  Then,
  \begin{equation*}
    \frac{1}{r_5} = \frac{\alpha - \frac{s_0}{\sigma}}{r_4} + \frac{\frac{s_0}{\sigma}}{r_6}.
  \end{equation*}
  We also have
  \begin{gather*}
    \frac{1}{r_{0,b}} = \frac{\alpha}{r_4} + \frac{1}{\rho_0} = \frac{1}{r_4} + \frac{1}{r_5}, \quad \frac{1}{r_{s_0,b}} = \frac{\alpha+1}{r_4} \\
    \frac{1}{\widetilde{\rho}_0^{\prime}} = \frac{1}{p_{0,b}} + \frac{1}{r_{0,b}} = \frac{1}{p_{s_0,b}} + \frac{1}{r_{s_0,b}}, \quad \frac{1}{\widetilde{\gamma}_0^{\prime}} = \frac{\alpha}{\gamma_1} + \frac{1}{\gamma_0}.
  \end{gather*}
  By these relations, the fractional Leibniz rule, \refsp{\lem:chain2}, and the Sobolev--Lorentz embedding, we obtain
  \begin{align*}
    & \norm{ D^{s_0} \left( a(x)|u|^{\alpha}w \right) }_{L^{\widetilde{\gamma}_0^{\prime}}(I,L^{\widetilde{\rho}_0^{\prime},2})} \\
    & \hspace{30pt} \lesssim \Big\| \norm{a(x)}_{L^{p_{0,b},\infty}} \norm{ D^{s_0} \left( |u|^{\alpha}w \right) }_{L^{r_{0,b},2}} \\
    & \hspace{60pt} + \norm{D^{s_0} a(x)}_{L^{p_{s_0,b},\infty}} \norm{ |u|^{\alpha}w }_{L^{r_{s_0,b},2}} \Big\|_{L^{\widetilde{\gamma}_0^{\prime}}(I)} \\
    & \hspace{30pt} \lesssim \Big\| \norm{u}_{L^{r_4,2}}^{\alpha} \norm{ D^{s_0} w }_{L^{\rho_0,2}} + \norm{ D^{s_0} (|u|^{\alpha})}_{L^{r_5,\frac{2}{\alpha}}} \norm{ w }_{L^{r_4,2}} \\
    & \hspace{180pt} + \norm{u}_{L^{r_4,2}}^{\alpha} \norm{w}_{L^{r_4,2}} \Big\|_{L^{\widetilde{\gamma}_0^{\prime}}(I)} \\
    & \hspace{30pt} \lesssim \norm{ \norm{u}_{L^{r_4,2}}^{\alpha - \frac{s_0}{\sigma}} \norm{D^{\sigma}u}_{L^{r_6,2}}^{\frac{s_0}{\sigma}} \norm{ D^{s_0} w }_{L^{\rho_0,2}} + \norm{u}_{H_{\rho_1,2}^{s}}^{\alpha} \norm{D^{s_0} w}_{L^{\rho_0,2}} }_{L^{\widetilde{\gamma}_0^{\prime}}(I)} \\
    & \hspace{30pt} \lesssim \norm{ \norm{u}_{H_{\rho_1,2}^{s}}^{\alpha} \norm{D^{s_0} w}_{L^{\rho_0,2}} }_{L^{\widetilde{\gamma}_0^{\prime}}(I)} \\
    & \hspace{30pt} \lesssim \norm{u}_{L^{\gamma_1}(I;H_{\rho_1,2}^{s})}^{\alpha} \norm{D^{s_0} w}_{L^{\gamma_0}(I,L^{\rho_0,2})}. 
  \end{align*}
  This proves \refsp{eq:5-2}.
  
  { The same estimate holds when $|u|^{\alpha}$ is replaced by an $\alpha$-H\"older continuous function $F(u)$ satisfying $F(0) = 0$. In particular, by the mean value formula, we obtain
  \begin{align*}
    & \norm{D^{s_0} \left( a(x)(|u|^{\alpha}u - |v|^{\alpha}v) \right)}_{L^{\widetilde{\gamma}_0^{\prime}}(I,L^{\widetilde{\rho}_0^{\prime},2})} \\
    & \hspace{30pt} \lesssim \left( \norm{u}_{L^{\gamma_1}(I,H_{\rho_1,2}^{s})} + \norm{v}_{L^{\gamma_1}(I,H_{\rho_1,2}^{s})} \right)^{\alpha} \norm{D^{s_0} (u - v)}_{L^{\gamma_0}(I,L^{\rho_0,2})}.
  \end{align*}
  }

  We next prove \refsp{eq:5-3}. Put 
  \begin{gather*}
    \frac{1}{r_7} \coloneq \frac{1}{\rho} - \frac{s}{N} = \frac{1}{\rho_2} - \frac{s_0}{N} \\
    \frac{1}{\widetilde{r}_{0,b}} \coloneq \frac{\alpha}{r_7} + \frac{1}{\rho}, \quad \frac{1}{\widetilde{r}_{s,b}} \coloneq \frac{\alpha + 1}{r_7}.
  \end{gather*}
  Note that 
  \begin{equation*}
    \frac{1}{r^{\prime}} = \frac{1}{p_{0,b}} + \frac{1}{\widetilde{r}_{0,b}} = \frac{1}{p_{s,b}} + \frac{1}{\widetilde{r}_{s,b}}.
  \end{equation*}

  { When $0<s<1$}, the fractional Leibniz rule, \refsp{\prop:Fujiwara}, and the Sobolev--Lorentz embedding give
  \begin{align*}
    & \norm{D^s\left[ a(x)\left( |u|^{\alpha}u - |v|^{\alpha}v \right) \right]}_{L^{r^\prime,2}} \\
    & \hspace{30pt} \lesssim \norm{D^s\left( |u|^{\alpha}u - |v|^{\alpha}v \right)}_{L^{\widetilde{r}_{0,b},2}} + \norm{|u|^{\alpha}u - |v|^{\alpha}v}_{L^{\widetilde{r}_{s,b},2}} \\
    & \hspace{30pt} \lesssim \left( \norm{u}_{L^{r_7,2}} + \norm{v}_{L^{r_7,2}} \right)^{\alpha} \norm{D^s(u-v)}_{L^{\rho,2}} \\
    & \hspace{60pt} + \norm{u-v}_{L^{r_7,2}}^{\alpha} \left( \norm{D^s u}_{L^{\rho,2}} + \norm{D^s v}_{L^{\rho,2}} \right) \\
    & \hspace{90pt} + \left( \norm{u}_{L^{r_7,2}} + \norm{v}_{L^{r_7,2}} \right)^{\alpha} \norm{u-v}_{L^{r_7,2}} \\
    & \hspace{30pt} \lesssim \left( \norm{u}_{H_{\rho,2}^{s}} + \norm{v}_{H_{\rho,2}^{s}} \right)^{\alpha} \norm{u-v}_{H_{\rho,2}^{s}} \\
    & \hspace{60pt} + \norm{D^{s_0}(u-v)}_{L^{\rho_2,2}}^{\alpha} \left( \norm{u}_{H_{\rho,2}^{s}} + \norm{v}_{H_{\rho,2}^{s}} \right).
  \end{align*}
  Applying H\"older's inequality in time yields \refsp{eq:5-3}. 

  When $s=1$, the same argument follows from the classical chain rule and the inequality
  \begin{equation*}
    |\nabla(|u|^{\alpha}u - |v|^{\alpha}v)| \lesssim (|u| + |v|)^{\alpha} |\nabla(u-v)| + |u-v|^{\alpha} (|\nabla u|+|\nabla v|).
  \end{equation*}

  For a general coefficient satisfying the \assumptionref, we apply the estimates above to the corresponding terms after the fractional Leibniz decomposition. The terms involving bounded components are treated as stated in \refsp{\rem:2-15}.

  This completes the proof of \refsp{\lem:5-3}.
\end{proof}

\begin{proof}[Proof of \refsp{\thm:1-2}]
  For the proof in the critical case, we first establish the following proposition.

  \begin{proposition} \label{proposition:5-4}
    Let $s\geq 0$, and let $\cK$ be a compact subset of $H^s$. Then, for any admissible pair $(q,r)$ satisfying $2<r<\infty$, 
    \begin{equation*}
      \lim_{T\to 0}\sup_{\varphi\in\cK}\norm{e^{i\cdot\Delta}\varphi}_{L^q(0,T;H_{r,2}^{s})} = 0.
    \end{equation*}
  \end{proposition}

  \begin{proof}
    By the Strichartz estimates \refsp{eq:14}, for any $\varphi\in\cK$, we have
    \begin{equation*}
      e^{it\Delta}\varphi\in L^q(\R;H_{r,2}^{s})
    \end{equation*}
    and since $q<\infty$, the absolute continuity of the integral implies
    \begin{equation*}
      \lim_{T\to 0}\norm{e^{i\cdot\Delta}\varphi}_{L^q(0,T;H_{r,2}^{s})} = 0.
    \end{equation*} 

    Let $\varepsilon>0$. For the constant $C$ in the Strichartz estimate, choose $\delta>0$ sufficiently small so that $C\delta<\frac{\varepsilon}{2}$. Since $\cK$ is compact, there exist $\varphi_1, \ldots, \varphi_k \in \cK$ such that $\cK\subset \displaystyle\bigcup_{j=1}^{k}B_{H^s}(\varphi_j, \delta)$. For each $j\in\{1,\ldots,k\}$, we can choose $\eta_j>0$ such that whenever $0<T<\eta_j$, 
    \begin{equation} \label{eq:5-4}
      \norm{e^{i\cdot\Delta}\varphi_j}_{L^q(0,T;H_{r,2}^{s})} < \varepsilon/2.
    \end{equation}
     Set $\eta = \displaystyle\min_{1\leq j\leq k}\{\eta_j\}$. Take $0<T<\eta$ and $\varphi\in\cK$. Then there is some $j_0\in \{1,\ldots, k\}$ such that $\varphi\in B_{H^s}(\varphi_{j_0}, \delta)$. Using \refsp{eq:5-4} and the Strichartz estimate, we obtain
    \begin{align*}
      \norm{e^{i\cdot\Delta}\varphi}_{L^q(0,T;H_{r,2}^{s})} &\leq \norm{e^{i\cdot\Delta}(\varphi - \varphi_{j_0})}_{L^q(0,T;H_{r,2}^{s})} + \norm{e^{i\cdot\Delta}\varphi_{j_0}}_{L^q(0,T;H_{r,2}^{s})} \\
      &\leq C\norm{\varphi - \varphi_{j_0}}_{H^s} + \frac{\varepsilon}{2} \\
      &\leq \varepsilon.
    \end{align*}
    The proof is complete.
  \end{proof}

  We first establish existence. With the later argument for continuous dependence in mind, we consider an arbitrary compact set $\cK\subset H^s$ and a function $u_0\in \cK$. Set $R_\cK\coloneq \sup_{\varphi\in \cK}\hsnorm{\varphi}$. When $b>0$, let $(\gamma^{(0)}, \rho^{(0)})$ be the admissible pair given by \refsp{\lem:4-1} with $b_{\ast} = 0$. We define $X(M_1,M_2,T)$ as the set of all functions $u$ satisfying
  \begin{equation*}
    u\in L^{\infty}(0,T;H^{s}) \cap L^{\gamma}(0,T;H_{\rho,2}^{s}) \cap L^{\gamma}(0,T;H_{\rho_{b_{\ast}},2}^{s}),
  \end{equation*}
  together with
  \begin{equation*}
    \norm{u}_{L_T^{\infty}H^{s}} + \norm{u}_{L_T^{\gamma^{(0)}}H_{\rho^{(0)},2}^{s}} \leq M_1, \quad \norm{u}_{L_T^{\gamma}H_{\rho,2}^{s}}\leq M_2.
  \end{equation*}
  For $u,v\in X(M_1,M_2,T)$, we further define the distance
  \begin{equation*}
    d(u,v) \coloneq \norm{u-v}_{L_T^{\infty}L^{2}} + \norm{u-v}_{L_T^{\gamma}L^{\rho,2}} + \norm{u-v}_{L_T^{\gamma^{(0)}}L^{\rho^{(0)},2}}.
  \end{equation*}
  When $b=0$, the terms involving $(\gamma^{(0)}, \rho^{(0)})$ in the definitions above are omitted. By the same argument as in \refsp{sec:4}, \((X(M_1,M_2,T),d)\) is a complete metric space.

  We first consider the terms corresponding to the exponent $b$. When $b>0$, the terms corresponding to $b_{\ast}=0$ are estimated by the subcritical argument in \refsp{sec:4}. These estimates contain a positive power of $T$, so they can be controlled by taking $T>0$ sufficiently small.

  We consider the mapping $\Phi$ defined in \refsp{eq:1-5}. By the Strichartz estimates \refsp{\prop:2} and the same nonlinear estimates as in the proof of \refsp{\lem:4-2}, with $\widetilde{s} = s$ for the critical component, we obtain
  \begin{align*}
    &\norm{\Phi(u)}_{L_T^{\infty}H^{s}} + \norm{\Phi(u)}_{L_T^{\gamma^{(0)}}H_{\rho^{(0)},2}^{s}} \\
    &\hspace{60pt} \leq C \left( \hsnorm{u_0} + \norm{u}_{L_T^{\gamma}H_{\rho,2}^{s}}^{\alpha+1} + T^{\beta_0}\norm{u}_{L_T^{\gamma^{(0)}}H_{\rho^{(0)},2}^{s}}^{\alpha + 1} \right),
  \end{align*}
  and
  \begin{equation*}
    \norm{\Phi(u)}_{L_T^{\gamma}H_{\rho,2}^{s}} \leq \norm{e^{it\Delta}u_0}_{L_T^{\gamma}H_{\rho,2}^{s}} + C \norm{u}_{L_T^{\gamma}H_{\rho,2}^{s}}^{\alpha+1} + CT^{\beta_0}\norm{u}_{L_T^{\gamma^{(0)}}H_{\rho^{(0)},2}^{s}}^{\alpha + 1},
  \end{equation*}
  where $C>0$ is a constant independent of $u_0$, $\cK$, $M_1$, $M_2$, and $T$.
  
  Choose $\eta>0$ sufficiently small so that $C(2\eta)^{\alpha}$ is small enough, and set $M_2\coloneq 2\eta$. We also set $M_1\coloneq 2C(R_\cK+M_2^{\alpha+1})$. By \refsp{\prop:5-4}, we may choose $T=T(\cK)>0$ sufficiently small so that
  \begin{equation*}
    \norm{e^{it\Delta}u_0}_{L_T^{\gamma}H_{\rho,2}^{s}} \leq \eta
  \end{equation*}
  for every $u_0\in \cK$. Decreasing $T$ further if necessary, we also assume that
  \begin{equation*}
    CT^{\beta_0}M_1^{\alpha + 1} \leq \min\left\{ \frac{M_1}{2}, \frac{\eta}{2} \right\}. 
  \end{equation*}
  Then, for any $u,v\in X(M_1,M_2,T)$ we obtain
  \begin{equation*}
    \norm{\Phi(u)}_{L_T^{\infty}H^{s}} + \norm{\Phi(u)}_{L_T^{\gamma^{(0)}}H_{\rho^{(0)},2}^{s}} \leq M_1, \quad \norm{\Phi(u)}_{L_T^{\gamma}H_{\rho,2}^{s}}\leq M_2.
  \end{equation*}
  The contraction estimate follows by the same argument as in the proof of \refsp{\thm:1-1}, with $\widetilde{s}=s$ for the critical terms. Thus,
  \begin{equation*}
    d(\Phi(u),\Phi(v)) \lesssim (M_2^{\alpha} T^{\beta_0} M_1^{\alpha}) d(u,v).
  \end{equation*}
  Since $M_2>0$ is sufficiently small, and $T>0$ can be taken smaller if necessary, the map $\Phi$ is a contraction on $X(M_1,M_2,T)$. The Banach fixed-point theorem guarantees the existence of a unique solution $u\in X(M_1,M_2,T)$. The time continuity of the Duhamel term also gives $u\in C([0,T];H^s)$. Finally, by shifting the initial time and repeating the same argument, we can extend the solution to its maximal existence time $T_{\max}(u_0)>0$.

  For (i), by the absolute continuity of the time integral, we can partition $[0,T]$ into finitely many subintervals $I_j$ such that
  \begin{equation*}
    \norm{u}_{L^{\gamma}(I_j;H_{\rho,2}^{s})} + \norm{v}_{L^{\gamma}(I_j;H_{\rho,2}^{s})}
  \end{equation*}
  is sufficiently small on each $I_j$. As in the contraction argument, the term involving $u-v$ on the right-hand side can be absorbed into the left-hand side. Repeating this procedure over the finite partition, we obtain uniqueness as in the subcritical case. For (ii), the same conclusion as in the subcritical case follows from the Strichartz estimates and the nonlinear estimates above.
  
  Next we consider (iii). For $\alpha \geq 1$, the assertion follows as in the subcritical case by using the smallness of the space-time norm obtained in the proof of existence. Therefore, we consider only the case $0<\alpha<1$. 

  When $b>0$, the terms corresponding to $b_{\ast}=0$ can be estimated as in \refsp{sec:4}, using the Sobolev--Lorentz embedding when necessary. These estimates contain a positive power of $T$, so we omit their details and focus on the terms corresponding to $b$.
    
  Fix $\nu>0$ sufficiently small. By the definition of $(\gamma_1,\rho_1)$ in \refsp{\lem:5-2} and the interpolation inequality \refsp{eq:2-3}, we have 
  \begin{equation*}
    \norm{u}_{L_T^{\gamma_1}H_{\rho_1,2}^{s}} \lesssim \norm{u}_{L_T^{\infty}H^{s}}^{1-\theta} \norm{u}_{L_T^{\gamma}H_{\rho,2}^{s}}^{\theta}.
  \end{equation*}
  Since $M_1$ remains bounded as $M_2\to 0$, we can choose $M_2>0$ sufficiently small and then $T=T(\cK)>0$ as in the existence argument so that
  \begin{equation*}
    \norm{u}_{L_T^{\gamma}H_{\rho,2}^{s}} + \norm{v}_{L_T^{\gamma}H_{\rho,2}^{s}} + \norm{u}_{L_T^{\gamma_1}H_{\rho_1,2}^{s}} + \norm{v}_{L_T^{\gamma_1}H_{\rho_1,2}^{s}} \leq \nu.
  \end{equation*}

  Arguing as in the contraction mapping argument, we first obtain
  \begin{align*}
    d(u,v) \lesssim \lnorm{u_0 - v_0} + (\nu^{\alpha} + T^{\beta_0}M_1^{\alpha}) d(u,v).
  \end{align*}
  Taking $\nu>0$ sufficiently small, the last term on the right-hand side can be absorbed into the left-hand side. Hence,
  \begin{equation}
    d(u,v) \lesssim \lnorm{u_0 - v_0}. \label{eq:5-5}
  \end{equation}
  Next, \refsp{\prop:2-12} and the difference estimate following \refsp{eq:5-2} yield
  \begin{align*}
    & \norm{D^{s_0}(u-v)}_{L_T^{\gamma_0}L^{\rho_0,2}} \\
    & \hspace{30pt} \lesssim \norm{D^{s_0}(u_0-v_0)}_{\dot{H}^{s - s_0}} + \nu^{\alpha} \norm{D^{s_0}(u-v)}_{L_T^{\gamma_0}L^{\rho_0,2}} \\
    & \hspace{30pt} = \lnorm{D^{s}(u_0-v_0)} + \nu^{\alpha} \norm{D^{s_0}(u-v)}_{L_T^{\gamma_0}L^{\rho_0,2}}.
  \end{align*} 
  Reducing $\nu>0$ further if needed, we obtain
  \begin{equation*}
    \norm{D^{s_0}(u-v)}_{L_T^{\gamma_0}L^{\rho_0,2}} \lesssim \hsnorm{u_0 - v_0}.
  \end{equation*}
  Furthermore, { choose an admissible pair $(q_{\ast},r_{\ast})\in\cA_0$ such that $r_{\ast} > \rho$ and $r_{\ast}$ is sufficiently close to $\rho$. Define $R_{\ast}$ by 
  \begin{equation*}
    \frac{1}{R_{\ast}} = \frac{1}{r_{\ast}} - \frac{s-s_0}{N}.
  \end{equation*}
  By the definitions of the exponents in \refsp{\lem:5-2}, there exists $0<\widetilde{\theta}<1$ such that 
  \begin{equation*}
    \frac{1}{\gamma_2} = \frac{\widetilde{\theta}}{\gamma_0} + \frac{1 - \widetilde{\theta}}{q_\ast}, \quad \frac{1}{\rho_2} = \frac{\widetilde{\theta}}{\rho_0} + \frac{1 - \widetilde{\theta}}{R_\ast}.
  \end{equation*}
  Thus, by interpolation and the Sobolev--Lorentz embedding,}
  \begin{align*}
    & \norm{D^{s_0}(u-v)}_{L_T^{\gamma_2}L^{\rho_2,2}} \\
    & { \hspace{50pt} \lesssim \norm{D^{s_0}(u-v)}_{L_T^{\gamma_0}L^{\rho_0,2}}^{\widetilde{\theta}} \norm{D^{s}(u-v)}_{L_T^{q_{\ast}}L^{r_{\ast},2}}^{1 - \widetilde{\theta}}} \\
    & \hspace{50pt} \lesssim_{\cK} \hsnorm{u_0 - v_0}^{\widetilde{\theta}}.
  \end{align*}
  Consequently, the Duhamel formula, the Strichartz estimates, and \refsp{eq:5-3} give
  \begin{align*}
    \hspace{30pt} &\norm{D^s(u - v)}_{L_T^{\infty}L^2} + \norm{D^s(u - v)}_{L_T^{\gamma}L^{\rho,2}} \\
    &\hspace{30pt} \lesssim \hsnorm{u_0 - v_0} + \nu^{\alpha} \norm{u-v}_{L_T^{\gamma}H_{\rho,2}^{s}} + \nu\norm{D^{s_0}(u-v)}_{L_T^{\gamma_2}L^{\rho_2,2}}^{\alpha} \\
    &\hspace{30pt} \lesssim_{\cK} \hsnorm{u_0 - v_0} + \nu^{\alpha} \norm{u-v}_{L_T^{\gamma}H_{\rho,2}^{s}} + \nu\hsnorm{u_0 - v_0}^{\widetilde{\theta}\alpha}.
  \end{align*}
  Combining \refsp{eq:5-5} with the above estimate, and absorbing the term multiplied by $\nu^{\alpha}$ into the left-hand side again, we obtain
  \begin{align*}
    \hspace{20pt} \norm{u - v}_{L_T^{\infty}H^s} + \norm{u  - v}_{L_T^{\gamma}H_{\rho,2}^{s}} \lesssim_{\cK} \hsnorm{u_0 - v_0} + \hsnorm{u_0 - v_0}^{\widetilde{\theta}\alpha}.
  \end{align*}
  Finally, applying the same argument to an arbitrary admissible pair $(\widetilde{q},\widetilde{r})\in\cA_0$, we obtain \refsp{eq:cri}, which proves the continuous dependence on the initial data.
\end{proof}

\section*{Acknowledgments}

The author would like to thank Professor Mishio Kawashita for his valuable guidance and helpful comments on this work. This work was also supported by JST SPRING, Grant Number JPMJSP2132.


\begin{thebibliography}{99}

\bibitem{A-T}
Aloui, L. and Tayachi, S.
Local well-posedness for the inhomogeneous nonlinear Schr\"o\-dinger equation.
\textit{Discrete Contin. Dyn. Syst.} 41 (11) (2021) 5409--5437.

\bibitem{B-D-S}
Bhimani, D. G., Dhingra, D., and Sohani, V. K.
Inhomogeneous nonlinear Schrö\-dinger equation in Fourier-Lebesgue and modulation spaces,
arXiv:2608.30777 [math.AP] (2026).

\bibitem{B-K-O}
Byeon J., Kim H. and Oh J., 
Interpolation inequalities in function spaces of Sobolev-Lorentz type, 
J. Math. Anal. Appl. {\bf 516} (2022), no.~2, Paper No. 126519, 37 pp.

\bibitem{Caz} 
T. Cazenave, 
{\it Semilinear Schr\"odinger equations}, 
Courant Lecture Notes in Mathematics, 10, New York Univ., Courant Inst. Math. Sci., New York, 2003 Amer. Math. Soc., Providence, RI, 2003

\bibitem{C-C-F}
Campos L., Correia S., and Farah L.~G., 
Sharp well-posedness and ill-posedness results for the inhomogeneous NLS equation, 
Nonlinear Anal. Real World Appl. {\bf 85} (2025), Paper No. 104336, 21 pp.

\bibitem{C-W}
Cazenave, T. and Weissler, F. B.,
Rapidly decaying solutions of the nonlinear Schr\"odinger equation, 
Comm. Math. Phys. {\bf 147} (1992), no.~1, 75--100.

\bibitem{C-N}
Cruz-Uribe, D. and Naibo, V.
Kato--Ponce inequalities on weighted and variable Lebesgue spaces.
\textit{Differential Integral Equations}, 29 (2016), 801--836.

\bibitem{D-F}
de~Bouard, A. and Fukuizumi, R., 
Stability of standing waves for nonlinear Schr\"odinger equations with inhomogeneous nonlinearities, 
Ann. Henri Poincar\'e{} {\bf 6} (2005), no.~6, 1157--1177.

\bibitem{F-W}
Fibich, G. and Wang, X.~P., 
Stability of solitary waves for nonlinear Schr\"odinger equations with inhomogeneous nonlinearities, 
Phys. D {\bf 175} (2003), no.~1-2, 96--108.

\bibitem{fos}
Foschi, D.,
Inhomogeneous Strichartz estimates, 
J. Hyperbolic Differ. Equ. {\bf 2} (2005), no.~1, 1--24.

\bibitem{fuj}
Fujiwara, K.,
Remark on the chain rule of fractional derivative in the Sobolev framework.
\textit{Math. Inequal. Appl.} 24 (4) (2021) 1113--1124.

\bibitem{G-S}
Genoud, F. and Stuart, C.~A., 
Schr\"odinger equations with a spatially decaying nonlinearity: existence and stability of standing waves, 
Discrete Contin. Dyn. Syst. {\bf 21} (2008), no.~1, 137--186.

\bibitem{gil} 
Gill, T.S., 
Optical guiding of laser beam in nonuniform plasma, 
Pramana J. Phys. {\bf 55} (5–6) (2000) 835–842.

\bibitem{G-V}
Ginibre, J., and Velo, G., 
The global Cauchy problem for the nonlinear Schr\"odinger equation revisited, 
Ann. Inst. H. Poincar\'e{} Anal. Non Lin\'eaire {\bf 2} (1985), no.~4, 309--327.

\bibitem{gra}
Grafakos, L.,
\textit{Classical Fourier analysis},
Springer, New York, 2008.

\bibitem{G-M-T}
Grira, M., Miyazaki, H. and Tayachi, S.,
Threshold for the existence of scattering states for inhomogeneous nonlinear Schr\"odinger equations without gauge invariance,
arXiv:2609.06404 [math.AP] (2026).

\bibitem{Guz} 
Guzm\'an, C.~M., 
On well posedness for the inhomogeneous nonlinear Schr\"odinger equation, 
Nonlinear Anal. Real World Appl. {\bf 37} (2017), 249--286

\bibitem{H-Y-Z}
Hajaiej H., Yu X. and Zhai Z., 
Fractional Gagliardo-Nirenberg and Hardy inequalities under Lorentz norms, 
J. Math. Anal. Appl. {\bf 396} (2012), no.~2, 569--577.

\bibitem{kat}
Kato, T.,
An $L^{q,r}$-theory for nonlinear Schr\"odinger equations, 
in {\it Spectral and scattering theory and applications}, 223--238, Adv. Stud. Pure Math., 23, Math. Soc. Japan, Tokyo.

\bibitem{K-T}
Keel, M. and Tao, T.~C., 
Endpoint Strichartz estimates, 
Amer. J. Math. {\bf 120} (1998), no.~5, 955--980.

\bibitem{K-L-S}
Kim, J., Lee, Y. and Seo, I., 
On well-posedness for the inhomogeneous nonlinear Schr\"odinger equation in the critical case, 
J. Differential Equations {\bf 280} (2021), 179--202.

\bibitem{lee}
Lee, Y., 
The 3D energy-critical inhomogeneous nonlinear Schr\"odinger equation with strong singularity, 
Nonlinear Anal. {\bf 272} (2026), Paper No. 114181, 18 pp.

\bibitem{L-Y-Y}
Liu, J., Yang, D.~C. and Yuan, W., 
Littlewood-Paley characterizations of anisotropic Hardy-Lorentz spaces, 
Acta Math. Sci. Ser. B (Engl. Ed.) {\bf 38} (2018), no.~1, 1--33.

\bibitem{one}
O'Neil, R.
Convolution operators and $L(p,q)$ spaces.
\textit{Duke Math. J.}, 30 (1963), 129--142.

\bibitem{par}
Park, D.,
Global well-posedness and scattering for the focusing energy-critical
inhomogeneous nonlinear Schrödinger equation with non-radial data,
arXiv:2410.12321 [math.AP] (2024).

\bibitem{str}
Strichartz, R.~S., 
Restrictions of Fourier transforms to quadratic surfaces and decay of solutions of wave equations, 
Duke Math. J. {\bf 44} (1977), no.~3, 705--714.

\bibitem{tag}
Taggart, R.~J., 
Inhomogeneous Strichartz estimates, 
Forum Math. {\bf 22} (2010), no.~5, 825--853.

\bibitem{vil}
Vilela, M.~C., 
Inhomogeneous Strichartz estimates for the Schr\"odinger equation, 
Trans. Amer. Math. Soc. {\bf 359} (2007), no.~5, 2123--2136.

\bibitem{yaj}
Yajima, K., 
Existence of solutions for Schr\"odinger evolution equations, 
Comm. Math. Phys. {\bf 110} (1987), no.~3, 415--426.

\end{thebibliography}
\end{document}